\documentclass{amsart}
\usepackage{amssymb}
\usepackage{euscript}
\usepackage[matrix,arrow,curve]{xy}
\usepackage{graphics}
\usepackage{color}
\usepackage{epsfig}

\renewcommand{\subsection}[1]{\hspace{-\parindent}\refstepcounter{subsection}{\bf
(\arabic{section}\alph{subsection}) #1.}}

\numberwithin{equation}{section}

\theoremstyle{plain}

\newtheorem*{theorem*}{Theorem}
\newtheorem{thm}{Theorem}[section]
\newtheorem{theorem}[thm]{Theorem}

\newtheorem{corollary}[thm]{Corollary}
\newtheorem{lemma}[thm]{Lemma}

\newtheorem{proposition}[thm]{Proposition}

\newtheorem{remark}[thm]{Remark}

\newcommand{\R}{\mathbb{R}}
\newcommand{\Z}{\mathbb{Z}}
\newcommand{\Q}{\mathbb{Q}}
\newcommand{\C}{\mathbb{C}}

\newcommand{\half}{{\textstyle\frac{1}{2}}}

\newcommand{\iso}{\cong}           %isomorphism sign
\newcommand{\CP}[1]{\C {\mathrm P}^{#1}}

\renewcommand{\O}{\EuScript{O}}
\newcommand{\g}{\mathfrak{g}}
\newcommand{\h}{\mathfrak{h}}

\newcommand{\Sym}{\mathrm{Sym}}
\newcommand{\Hom}{\mathrm{Hom}}
\newcommand{\id}{\mathrm{id}}
\newcommand{\F}{\EuScript{F}}
\newcommand{\A}{\EuScript{A}}
\newcommand{\B}{\EuScript{B}}

\newcommand{\scrM}{\EuScript{M}}

\newcommand{\scrL}{\EuScript{L}}

\begin{document}
\title{Homological mirror symmetry for the genus two curve}
\author{Paul Seidel}
%\date{v6; October 8, 2009}
\maketitle

\section{Introduction}

Homological Mirror Symmetry (HMS) relates algebraic and symplectic geometry through their associated categorical structures. This relation is by no means straightforward, and its exploration has been the driving force behind several recent developments in algebraic geometry. Let's consider the simplest case, that of an elliptic curve. There, HMS provides a new viewpoint \cite{subotic} on the well-known classification of vector bundles and their tensor category structure \cite{atiyah57}. It also gives rise to nontrivial identities involving theta series for indefinite quadratic forms \cite{polishchuk98}. Finally, it sheds light on the action of $\mathit{SL}(2,\Z)$ on the derived category of an elliptic curve, first introduced in \cite{mukai}. Higher-dimensional instances of HMS also lead to consequences of the same kind, even though those tend to be somewhat harder. For instance, many symplectic manifolds have rich automorphism groups arising from monodromy, and by thinking of those in terms of HMS, new classes of autoequivalences of derived categories have been discovered \cite{seidel-thomas, horja, huybrechts-thomas}. A more recent example is the discovery of a new relation between coherent and constructible sheaves, in the context of toric varieties \cite{4}. Other areas that have been influenced by HMS include: the study of spaces of stability conditions; tropical geometry; the theory of exceptional collections and mutation; and singularity theory.

Having considered its context, we'd now like to concentrate on the HMS conjecture itself (the following discussion is definitely not intended to be complete: it only lists a few papers which are relevant for the developments presented in the body of the paper). The conjecture comes in several related but distinct versions, which apply in different geometric contexts. Kontsevich's original version \cite{kontsevich94} concerned Calabi-Yau varieties. There, we now have complete proofs of some instances \cite{polishchuk-zaslow98,seidel03b} (including the elliptic curve case mentioned above) and partial results for many more \cite{kontsevich-soibelman00,fukaya02b}. Soon after, Kontsevich himself proposed an analogous conjecture for Fano varieties. This was gradually extended further, and it seems that varieties with effective anticanonical divisor provide a natural context \cite{auroux07}. The mirror in this case is not another variety but rather a Landau-Ginzburg theory, which means a variety together with a holomorphic function. Because of this asymmetry, the two directions of the mirror correspondence lead to substantially different mathematics. The one relevant for our purpose is where the Landau-Ginzburg theory is considered algebro-geometrically, through matrix factorizations or more generally Orlov's Landau-Ginzburg branes \cite{orlov04}. The corresponding symplectic geometry has been addressed in the toric case in \cite{cho02, cho-oh02, fooo08,fooo08b,auroux07}. To the best of this author's knowledge, none of those papers actually proves HMS in its full form, but in many cases it should follow from the results presented there together with additional steps which are fundamentally well-understood.

More recently, Katzarkov \cite{katzarkov07, katzarkov-kontsevich-pantev08} has proposed a further extension of HMS encompassing some varieties of general type. As before, the mirror is a Landau-Ginzburg theory. Abouzaid, Auroux, Gross, Katzarkov, and Orlov have explored both directions of the correspondence, and accumulated large amounts of evidence ($K$-theory computations \cite{abouzaid08,orlov08} and more unpublished material). The aim of this paper is to prove one direction of Katzarkov's conjecture in the simplest possible case. This is inspired by the work we've just mentioned, and additionally by another instance of mirror symmetry in the literature, namely genus zero curves with three orbifold points (see \cite[Section 7]{takahashi07}, and a little more recently \cite{rossi08}). Let $M$ be a genus two curve, equipped with a symplectic structure. Its mirror is a three-dimensional Landau-Ginzburg theory $X \rightarrow \C$, whose zero fibre $H \subset X$ is the union of three rational surfaces. The singular set $\mathrm{Sing}(H)$ is the union of three rational curves, intersecting as shown in Figure \ref{fig:curves} (the fact that this ``looks like'' a degeneration of the genus two curve is explained by an unpublished result of Gross-Katzarkov, which identifies the cohomology of $M$ with that of the sheaf of vanishing cycles on $H$). Details of the construction of the mirror will be given later.
\begin{figure}
\begin{centering}
\begin{picture}(0,0)%
\includegraphics{theta.pstex}%
\end{picture}%
\setlength{\unitlength}{3947sp}%
\begingroup\makeatletter\ifx\SetFigFont\undefined%
\gdef\SetFigFont#1#2#3#4#5{%
  \reset@font\fontsize{#1}{#2pt}%
  \fontfamily{#3}\fontseries{#4}\fontshape{#5}%
  \selectfont}%
\fi\endgroup%
\begin{picture}(622,1524)(1938,-973)
\end{picture}%
\caption{\label{fig:curves}}
\end{centering}
\end{figure}

Let $\F(M)$ be the Fukaya category of $M$, and $D^\pi(\F(M))$ its split-closed (Karoubi completed) derived category. On the other side, take $D^b_{\mathrm{sing}}(H)$ to be the category of Landau-Ginzburg branes, and let $D^\pi_{\mathrm{sing}}(H)$ be the split-closure of that.

\begin{theorem} \label{th:main}
There is an equivalence of triangulated categories,
\begin{equation}
D^\pi(\F(M)) \iso D^\pi_{\mathrm{sing}}(H).
\end{equation}
\end{theorem}

As one consequence, we have the maybe surprising result that the genus two mapping class group acts faithfully on $D^\pi_{\mathrm{sing}}(H)$. An outline of the proof of Theorem \ref{th:main} is given in the next section. For now, it is maybe enough to say that the argument relies on the fact that both categories can be described by $A_\infty$-algebras of a very special form ($A_\infty$-deformations of the exterior algebra, with an added group action). The determination of the exact $A_\infty$-structure is then reduced to the computation of finitely terms of a superpotential. In a wider context, Theorem \ref{th:main} may raise more questions than it answers. First of all, Katzarkov's original construction, which embeds $M$ as a holomorphic curve into $\CP{1} \times \CP{1}$ and applies \cite{hori-vafa00} to that situation, leads to a mirror which is similar but not quite the same as the one considered here. Presumably, the two resulting categories $D^\pi_{\mathrm{sing}}$ are equivalent, but that remains to be shown. Next, the approach followed here has a natural generalization to higher genus curves, and to some higher-dimensional manifolds. In a different direction, the use of split-closures is unsatisfactory, since that process is known to lose information \cite{orlov08b}.

{\em Addenda:} Several relevant preprints have appeared in the time since this one was originally written and distributed. A general approach to matrix factorizations similar to that in Sections \ref{sec:koszul}--\ref{sec:matrix} is given in \cite{dyckerhoff09}. Particularly relevant for us is \cite[Theorem 4.3]{dyckerhoff09}, which could replace the ad hoc computation \cite{seidel08x} in the proof of Proposition \ref{th:endomorphism-characterized}. The results of \cite{mehrotra05} and \cite{quintero-velez08} are generalized and put on a more systematic footing in \cite{baranovsky-pecharich09}. In particular, \cite[Theorem 2.10]{baranovsky-pecharich09} covers our needs in that respect, since it directly implies Theorem \ref{th:derived-mckay}. The generalization of our results to curves of any genus $\geq 2$ is carried out in \cite{efimov09}. Finally, \cite{kapustin-katzarkov-orlov-yotov09} contains, among other things, a detailed description of the mirror geometry in Katzarkov's original construction.

{\em Acknowledgments:} The author would like to thank Mohammed Abouzaid, Denis Auroux, Ludmil Katzarkov, Dima Orlov, and Ivan Smith for their help, as well as the NSF for partial support through grant DMS-0652620.

\section{Overview}
We now give a guided tour of the proof, simultaneously fixing the notation. Take $V = \C^3$. We write $\xi_k$ for the standard basis vectors of $V$, thought of as constant vector fields, and $v_k \in V^\vee$ for the dual basis of functions. The superpotential which is key to our considerations is the polynomial
\begin{equation} \label{eq:w}
W = - v_1v_2v_3 + v_1^5 + v_2^5 + v_3^5 \in \C[V].
\end{equation}
Take $Z \iso \Z/5$ to be the subgroup of $\mathit{SL}(V)$ generated by the diagonal matrix $\mathrm{diag}(\zeta,\zeta,\zeta^3)$, with $\zeta = \exp(2\pi i/5)$. Let $X \rightarrow \bar{X} = V/Z$ be the crepant resolution given by the $Z$-Hilbert scheme \cite{nakamura99}. Our mirror Landau-Ginzburg model is the composition
\begin{equation}
X \longrightarrow \bar{X} \xrightarrow{W} \C.
\end{equation}
In particular, $H \subset X$ is the preimage of $\bar{H} = W^{-1}(0)/Z \subset \bar{X}$. It is an elementary exercise to determine the geometry of $H$, which is as described in the Introduction. This is done in Section \ref{sec:mckay}.

At the same time, this construction yields a way to approach $D^\pi_{\mathrm{sing}}(H)$. Namely, a version of the derived McKay correpondence  \cite{mehrotra05, quintero-velez08} shows that this is equivalent to the equivariant category $D^\pi_{\mathrm{sing},Z}(W^{-1}(0))$. For simplicity, let's forget about the group action and just talk about $D^\pi_{\mathrm{sing}}(W^{-1}(0))$. A theorem of Orlov \cite{orlov08b} shows that this category is split-generated by a single object, which is the skyscraper sheaf at the origin, denoted by $\EuScript{S}_{W^{-1}(0),0}$. Hence, the category can be completely reconstructed from the $A_\infty$-structure on
\begin{equation}
\mathrm{Hom}_{D^\pi_{\mathrm{sing}}(W^{-1}(0))}(\EuScript{S}_{W^{-1}(0),0},\EuScript{S}_{W^{-1}(0),0}) \iso \Lambda(V).
\end{equation}
Moreover, matrix factorizations give rise to a natural dg structure underlying this algebra. On general grounds, the $A_\infty$-structure can be extracted from this by applying the Homological Perturbation Lemma (even though in practice, the computational complexity of computing the operations $\mu^d$ rises very rapidly with $d$). This is explained in Sections \ref{sec:koszul}--\ref{sec:matrix}.

Switching to the other side, we represent $M$ as a covering of a genus zero orbifold $\bar{M}$, where the covering group is $\Sigma = \mathrm{Hom}(Z,\C^*) \iso \Z/5$. We choose a collection of five curves $\{L_1,\dots,L_5\}$ which split-generate $D^\pi(\F(M))$, and which all project to the same immersed curve $\bar{L} \subset \bar{M}$. If we again forget about the covering group action, all the desired information is contained in the $A_\infty$-structure on the Floer cohomology
\begin{equation}
\mathit{HF}^*(\bar{L},\bar{L}) \iso \Lambda(V).
\end{equation}
The first few $A_\infty$-operations can be determined combinatorially by counting polygons (an idea that goes back at least to \cite{kontsevich94}). The relevant material is covered in Sections \ref{sec:fukaya}--\ref{sec:g2}.

At this point, we've reduced both sides to the computation of a specific $\Z/2$-graded $A_\infty$-deformation of the exterior algebra $\Lambda(V)$. The relevant deformation theory is governed by the differential graded Lie algebra of Hochschild cochains. We apply a version of Kontsevich's Formality Theorem \cite{kontsevich97}, and standard tools from Maurer-Cartan theory, to reduce this to a problem about polyvector fields, which means elements of $\C[[V]] \otimes \Lambda(V)$. In fact, in our case the $A_\infty$-deformation is determined by a single function $W \in \C[[V]]$, which turns out to be precisely the polynomial defined above. The crucial ingredient is a technical result (Proposition \ref{th:classification}) which shows that a specific isomorphism class of $A_\infty$-deformations, denoted by $\A$, is characterized by the first few nontrivial $A_\infty$-products. The underlying geometric idea is finite determinacy of function germs, which applies to any formal power series with an isolated critical point at the origin. This is the content of Sections \ref{sec:kontsevich}--\ref{sec:classify}.

\section{Kontsevich formality\label{sec:kontsevich}}

We begin by recalling some well-known generalities. Let $\g$ be a dg Lie algebra over $\C$. A Maurer-Cartan element is an $\alpha \in \g^1$ which satisfies
\begin{equation} \label{eq:mc0}
\partial\alpha + \half [\alpha,\alpha] = 0.
\end{equation}
There is a natural Lie algebra homomorphism from $\g^0$ to the space of affine vector fields on $\g^1$, which associates to $\gamma \in \g^0$ the infinitesimal gauge transformation $\alpha \mapsto -\partial\gamma  + [\gamma,\alpha]$. These endomorphisms are tangent to \eqref{eq:mc0}. Hence, in situations where they can be exponentiated, we get a group action on the space of solutions to the Maurer-Cartan equation. Two extreme special cases are worth considering. First, if $[\cdot,\cdot] = 0$, the Maurer-Cartan equation is just the cocycle equation, and infinitesimal gauge transformations act by adding coboundaries. On the other hand, if $\partial = 0$, what we have is just the adjoint action of $\g^0$ on the set of elements of $\g^1$ satisfying $[\alpha,\alpha] = 0$.

It is convenient to introduce generalized morphisms between dg Lie algebras, technically known as $L_\infty$-homo\-morphisms. Such a morphism $\Phi: \g \rightarrow \h$ consists of a sequence of multilinear maps $\Phi^k: \g^{\otimes k} \rightarrow \h$ of degree $1-k$, $k \geq 1$, which are antisymmetric in a suitably graded sense, and satisfy the equations spelled out in \cite{lada-markl95}. In particular, $\Phi^1$ is a chain map, and induces a homomorphism of graded Lie algebras on the cohomology level. One advantage is that quasi-isomorphisms can be inverted in this context. The precise statement we need is this:

\begin{lemma} \label{th:partial-reverse}
Let $\g$ be a graded Lie algebra, $\h$ a dg Lie algebra, and $\Psi: \g \rightarrow \h$ an $L_\infty$-homomorphism. Suppose that we are given a chain map $\Phi^1: \h \rightarrow \g$ and an $l: \h \rightarrow \h$ of degree $-1$, such that
\begin{equation}
\begin{aligned}
& \Phi^1 \circ \Psi^1 = \id, \\
& \Psi^1 \circ \Phi^1 - \id = \partial l + l \partial.
\end{aligned}
\end{equation}
Then $\Phi^1$ can be extended to an $L_\infty$-homomorphism $\Phi: \h \rightarrow \g$. Moreover, the higher order terms of $\Phi$ are given by universal formulae, which depend only on $\Psi$, $\Phi^1$ and $l$.
\end{lemma}

\proof The first part of the construction uses only $\Psi^1$, $\Phi^1$ and $l$. Given these, the Homological Perturbation Lemma constructs another $L_\infty$-algebra $\tilde{\g}$ with vanishing differential, whose underlying vector space is $\g$, together with an $L_\infty$-homomorphism $\tilde\Phi: \h \rightarrow \tilde{\g}$ whose first order term is $\Phi^1$. For explicit formulae see \cite{markl04} (that reference concerns $A_\infty$-algebras, but the $L_\infty$-case is parallel). Now $\tilde\Phi \circ \Psi: \g \rightarrow \tilde{\g}$ is an $L_\infty$-homomorphism whose first order term is the identity. Hence it admits a unique exact inverse with the same property. Define $\Phi = (\tilde\Phi \circ \Psi)^{-1} \circ \tilde\Phi$. \qed

To make the various pieces come together, we need to place ourselves in a framework where certain convergence properties are guaranteed. Assume that $\g$ is filtered pronilpotent, which means that it comes with a complete decreasing filtration $L_\bullet\g$ such that $L_1\g = \g$ and
\begin{equation}
\begin{aligned}
& \partial(L_r\g) \subset L_r\g,
& [L_r\g,L_s\g] \subset L_{r+s}\g.
\end{aligned}
\end{equation}
Then $\g^0$ is a pronilpotent Lie algebra, hence can be exponentiated to a prounipotent group by using the Baker-Campbell-Hausdorff formula. This group will act on the set of Maurer-Cartan elements. We call two elements equivalent if they lie in the same orbit (one can also define this relation through a suitable notion of homotopy between Maurer-Cartan elements). Deligne's basic idea (see \cite{goldman-millson88} and references therein) is that this is a good model for many kinds of deformation theory.

Let $\g,\h$ be two filtered pronilpotent dg Lie algebras. A filtered $L_\infty$-homo\-morphism $\Phi: \h \rightarrow \g$ consists of a family of maps as before, with the additional condition that
\begin{equation} \label{eq:phi-filtered}
\Phi^k(L_{r_k}\g \otimes \cdots \otimes L_{r_1}\g) \subset L_{r_k+\cdots+r_1}\g.
\end{equation}
There is an induced map on Maurer-Cartan elements, which preserves equivalence, namely:
\begin{equation} \label{eq:map-alpha}
\alpha \longmapsto \sum_{k=1}^\infty \textstyle\frac{1}{k!} \Phi^k(\alpha,\cdots,\alpha).
\end{equation}

\begin{lemma} \label{th:filtered}
Suppose that $\Phi: \h \rightarrow \g$ is a filtered $L_\infty$-homomorphism which is also a filtered quasi-isomorphism; the latter property means that $\Phi^1$ induces quasi-isomorphisms of chain complexes $L_r\h/L_{r+1}\h \rightarrow L_r\g/L_{r+1}\g$ for any $r$. Then \eqref{eq:map-alpha} induces a bijection between equivalence classes of Maurer-Cartan solutions.
\end{lemma}

This is an adapted version of a result from \cite[Section 4.4]{kontsevich97}. The proof given there does not immediately carry over to the filtered context. However, one can prove the result by a more direct obstruction theory computation, where solutions of the Maurer-Cartan equation (or homotopies between them) are lifted from $\g$ to $\h$ up to errors which are of successively higher order in our filtration.

Kontsevich \cite{kontsevich97} used this framework to explore the relation between commutative and noncommutative geometry. We will summarize his result, with some minor modifications. On the commutative geometry side, let's temporarily generalize our notation to allow $V = \C^n$ for any $n$. By definition, the space of formal polyvector fields on $V$ is
\begin{equation} \label{eq:poly}
\C[[V]] \otimes \Lambda(V) = \prod_{i,j} \mathrm{Sym}^i(V^\vee) \otimes \Lambda^j(V).
\end{equation}
The $(i,j)$ piece is given degree $j-1$, and the whole space becomes a graded Lie algebra with the Schouten bracket
\begin{equation}
\begin{aligned} {}
& [f \,\xi_{i_1} \wedge \cdots \wedge \xi_{i_k}, g \,\xi_{j_1} \wedge \cdots \xi_{j_l}] = \\
& \textstyle \qquad \sum_q (-1)^{k-q-1} f \, (\partial_{i_q}g) \,
\xi_{i_1}  \wedge \cdots \wedge \widehat{\xi_{i_q}} \wedge \cdots
\wedge \xi_{i_k} \wedge \xi_{j_1} \wedge \cdots \wedge \xi_{j_l} + \\ & \textstyle \qquad
 \sum_q (-1)^{l-q+ (k-1)(l-1)} g \,(\partial_{j_q}f)\, \xi_{j_1} \wedge \cdots \wedge \widehat{\xi_{j_q}} \wedge \cdots \wedge \xi_{j_l} \wedge \xi_{i_1} \wedge
 \cdots \wedge \xi_{i_k}.
\end{aligned}
\end{equation}
The Maurer-Cartan equation for $\alpha \in \C[[V]] \otimes \Lambda^2(V)$ says that the associated bracket $\{f,g\} = \alpha(df \wedge dg)$ satisfies the Jacobi identity, hence gives rise to a formal Poisson structure. Elements $\gamma \in \C[[V]] \otimes V$ are formal vector fields, acting by their Lie derivative. At least for vector fields vanishing at zero, the action can be exponentiated, and results in the obvious pushforward action of formal diffeomorphisms on Poisson brackets.

The analogue of polyvector fields in noncommutative geometry is given by Hochschild cohomology, which we now describe. Let $A$ be a graded associative algebra over $\C$. Its Hochschild complex $\mathit{CC}(A,A)$ is the space of graded multilinear maps
\begin{equation} \label{eq:cc-definition}
\mathit{CC}^d(A,A) = \prod_{i+j-1 =d} \mathrm{Hom}^j(A^{\otimes i},A).
\end{equation}
The Hochschild differential is
\begin{equation}
\begin{aligned}
& (\partial\phi)^j(a_j,\dots,a_1) = \textstyle\sum_k (-1)^{|\phi| + |a_1| + \cdots + |a_k| + k} \phi^{j-1}(a_j,\dots,a_{k+1}a_k,\dots,a_1) \\ & \qquad \qquad \qquad \qquad \quad
+ (-1)^{|\phi| + |a_1|+\cdots+|a_{j-1}| + j} a_j \phi^{j-1}(a_{j-1},\dots,a_1) \\ &
\qquad\qquad\qquad\qquad\quad
+ (-1)^{(|\phi|-1)(|a_1|-1)+1} \phi^{j-1}(a_j,\dots,a_2)a_1, \\
\end{aligned}
\end{equation}
and the Gerstenhaber bracket is
\begin{equation}
\begin{aligned}{}
[\phi,\psi]^j(a_j,\dots,a_1) & = \textstyle\sum_{k,l}
(-1)^{|\psi|(|a_1|+\cdots+|a_k|-k)} \phi^{j-l+1}(a_j,\dots, \\ & \qquad \qquad \qquad \qquad \quad \psi^l(a_{k+l},\dots,a_{k+1}),a_k,\dots,a_1) \\
&
- \textstyle\sum_{k,l}
(-1)^{|\phi| \cdot |\psi| + |\phi|(|a_1|+\cdots+|a_k|-k)} \psi^{j-l+1}(a_j,\dots, \\ & \qquad \qquad \qquad \qquad \quad \phi^l(a_{k+1},\dots,a_{k+1}),a_k,\dots,a_1).
\end{aligned}
\end{equation}
The cohomology of $\partial$ is just the Hochschild cohomology $\mathit{HH}(A,A)$, with the grading shifted down by $1$ from the standard convention. Take $\alpha \in \mathit{CC}^1(A,A)$, which by definition is a sequence of maps $\alpha^j: A^{\otimes j} \rightarrow A$ of degree $2-j$, $j \geq 0$. Set
\begin{equation}
\begin{cases} \mu^j = \alpha^j \text{ for $j \neq 2$}, \\
\mu^2(a_2,a_1) = \alpha^2(a_2,a_1) + (-1)^{|a_1|} a_2a_1.
\end{cases}
\end{equation}
Then, the Maurer-Cartan equation for $\alpha$ says that $\mu$ satisfies the equations for a curved $A_\infty$-structure, see for instance \cite{fooo}. Suppose for technical simplicity that $A$ is finite-dimensional in each degree, and take some $\gamma \in \mathit{CC}^0(A,A)$ whose constant term $\gamma^0 \in A^1$ vanishes. Define
\begin{equation} \label{eq:exp-gamma}
\begin{cases}
\phi^1 = \mathrm{id} + \gamma^1 + \half \gamma^1\gamma^1 + \cdots = \exp(\gamma^1), \\
\phi^2 = \gamma^2 + \half \gamma^1\gamma^2 + \half \gamma^2(\gamma^1 \otimes \mathrm{id}) +
\half \gamma^2(\mathrm{id} \otimes \gamma^1) + \textstyle\frac{1}{3} \gamma^2(\gamma^1 \otimes \gamma^1) + \cdots, \\
\dots
\end{cases}
\end{equation}
The general rule for $\phi^j$ is to sum up all possible ways of concatenating components of $\gamma$ to get a $j$-linear map. If there are $r$ components, and $s$ ways of ordering the components compatibly with their appearance in the concatenation, then the constant in front of the associated term is $s(r!)^{-1}$ (this in particular ensures convergence of the sums). If $\alpha$ and $\tilde\alpha$ are two Maurer-Cartan elements which are related by the exponentiated action of $\gamma$, the associated curved $A_\infty$-structures $\mu$, $\tilde\mu$ are related by $\phi$, which is an $A_\infty$-isomorphism.

We now specialize to exterior algebras $A = \Lambda(V)$. A classical result \cite{hochschild-kostant-rosenberg62} is that $\mathit{HH}(A,A) \iso \C[[V]] \otimes \Lambda(V)$. This isomorphism is induced by the Hochschild-Kostant-Rosenberg map, which is the projection $\Phi^1: \mathit{CC}(A,A) \rightarrow \mathrm{Hom}(T(V),\Lambda(V)) \rightarrow \C[[V]] \otimes \Lambda(V)$. Explicitly, thinking of $\Phi^1(\beta)$ as a $\Lambda(V)$-valued formal power series, we have
\begin{equation} \label{eq:hkr}
\Phi^1(\beta)(\xi) = \sum_{j=1}^{\infty} \beta^j(\xi,\dots,\xi).
\end{equation}
Kontsevich's formality theorem \cite{kontsevich97} says the following:

\begin{theorem} \label{th:kontsevich}
$\Phi^1$ is the first term of an $L_\infty$-homomorphism $\Phi$. Moreover, $\Phi$ is equivariant with respect to the action of $\mathit{GL}(V)$ on both sides.
\end{theorem}

Our formulation differs from the original one in two respects. First, it concerns exterior algebras instead of polynomial ones. However, the proof adapts in a straightforward way, by exchanging odd and even variables. Secondly, Kontsevich actually constructs an $L_\infty$-homomorphism $\Psi$ in the opposite direction. We want to use Lemma \ref{th:partial-reverse} to reverse direction, and that requires a choice of homotopy. By thinking of the classical grading of Hochschild cohomology, one sees that the homotopy can be taken to be a collection of maps $l^{i,j}: \mathrm{Hom}^j(A^{\otimes i},A) \longrightarrow \mathrm{Hom}^j(A^{\otimes i-1},A)$. Since each of these spaces is finite-dimensional and the group $\mathit{GL}(V)$ is reductive, one can average the homotopy to make it $\mathit{GL}(V)$-equivariant as well. Kontsevich's $\Psi$ is $\mathit{GL}(V)$-equivariant, and because of the way in which the inverse is defined, the same will then hold for $\Phi$.

Neither $\Lambda(V) \otimes \C[[V]]$ nor $\mathit{CC}(A,A)$ are pronilpotent, but one can remedy that by introducing an additional formal parameter $\hbar$, as in Kontsevich's original application to deformation quantization. This will also be the case here, but our parameter will have nonzero degree, which makes the situation somewhat less standard.

\section{Finite determinacy}

From now on, we again restrict to $V = \C^3$. Take the group $G \subset SL(V)$ which consists of diagonal matrices whose nonzero coefficients are fifth roots of unity. Because of the condition on the determinant, $G \iso (\Z/5)^2$. We will now tweak the previous framework by introducing the abovementioned formal parameter with nonzero degree, while at the same time adding equivariance with respect to $G$. Namely, let $\g$ be the graded vector space defined by
\begin{equation} \label{eq:g-algebra}
\g^d \;\; = \!\!\!\!\! \prod_{\substack{2i+j-4k = 3d+3 \\ k \geq 0, \, i \geq d+2}}\!\!\!\!\! (\Sym^i(V^\vee) \otimes \Lambda^j(V))^G\, \hbar^k.
\end{equation}
Even though the degrees are now different, their parities are the same as in our original discussion of $\C[[V]] \otimes \Lambda(V)$. The Nijenhuis bracket turns $\g$ into a graded Lie algebra. It is filtered pronilpotent, with $L_r\g^d$ being the part of \eqref{eq:g-algebra} where $i \geq d+1+r$.

For elements of $\g$ of a fixed degree $d$, the power of $\hbar$ occurring in each term $\Sym^i(V^\vee) \otimes \Lambda^j(V)$ is fixed to be $\frac{1}{4}(2i+j-3d-3)$. Hence, we can usually omit it from the notation, as long as we still remember the inequality
\begin{equation}
2i+j \geq 3d+3
\end{equation}
as well as the congruence
\begin{equation} \label{eq:mod4}
2i+j+d+1 \equiv 0\; \mathrm{mod}\; 4.
\end{equation}
Let $F_\bullet\C[[V]]$ be the complete decreasing filtration such that $F_r\C[[V]]$ consists of those power series with no terms of order strictly less than $r$. An element $\alpha \in \g^1$ has the form $\alpha = (\alpha^0,\alpha^2)$, where $\alpha^0 \in F_3\C[[V]]$ is an odd formal function, and $\alpha^2 \in F_4\C[[V]] \otimes \Lambda^2(V)$ is an even formal two-form. Here the terms even and odd refer to the action of $-1 \in \mathit{GL}(V)$ on polyvector fields. This property of $\alpha$ is an obvious consequence of \eqref{eq:mod4}. Similarly, an element $\gamma \in \g^0$ can be written as $\gamma = (\gamma^1,\gamma^3)$, where $\gamma^1 \in F_3\C[[V]] \otimes V$ is even, and $\gamma^3 \in F_2\C[[V]] \otimes \Lambda^3(V)$ is odd. On top of that, we of course have the $G$-invariance condition.

Before starting actual computations, it is worth while to acquire some geometric intuition. The Maurer-Cartan equation decomposes into
\begin{equation}
\half [\alpha^2,\alpha^2] = 0, \quad [\alpha^0,\alpha^2] = 0.
\end{equation}
As before, the first part says that $\alpha^2$ defines a Poisson bracket $\{\cdot,\cdot\}$. The second part says that $\{\alpha^0,\cdot\}$ is trivial, which means that the Poisson vector field associated to the function $\alpha^0$ is identically zero. Equivalently, $\alpha^2$ is a cocycle in the Koszul complex given by contraction with $d\alpha^0 \in \C[[V]] \otimes V^\vee$, which is
\begin{equation} \label{eq:koszul-complex}
0 \rightarrow \C[[V]] \otimes \Lambda^3(V) \xrightarrow{\iota_{d\alpha^0}} \C[[V]] \otimes \Lambda^2(V)
\xrightarrow{\iota_{d\alpha^0}} \C[[V]] \otimes V \xrightarrow{\iota_{d\alpha^0}} \C[[V]] \rightarrow 0.
\end{equation}
Considering degree zero elements, the exponentiated adjoint action of $\gamma = (\gamma^1,0)$ is the usual action of formal diffeomorphisms on polyvector fields. The exponentiated adjoint action of $\gamma = (0,\gamma^3)$ is given by
\begin{equation} \label{eq:gamma3-action}
(\alpha^0,\alpha^2) \longmapsto (\alpha^0,\alpha^2 + \iota_{d\alpha^0}\gamma^3).
\end{equation}

Note that $(W,0) \in \g^1$ is a solution of the Maurer-Cartan equation. It turns out that any other solution which is sufficiently close to this one is actually equivalent to it. The precise statement is:

\begin{lemma} \label{th:singularity}
Any Maurer-Cartan element $\alpha = (\alpha^0,\alpha^2) \in \g^1$ such that $\alpha^0 \equiv W$ mod $F_7\C[[V]]$ is equivalent to $(W,0)$.
\end{lemma}

Note that $W$ has an isolated singularity at the origin, which in algebraic terms means that the ideal $I = (\partial_1W,\dots,\partial_3 W) \subset \C[[V]]$ is of finite codimension. As a consequence, any formal power series which agrees with $W$ to sufficiently high order can be transformed into $W$ by a formal change of coordinates. This phenomenon is known in singularity theory as finite determinacy \cite{tougeron70} (see \cite[vol.\ I p.\ 121]{arnold-gusein-zade-varchenko} for an exposition). We will not appeal to these general results, but they've definitely guided our approach.

The explicit computation goes as follows. Elementary manipulation shows that
\begin{equation} \label{eq:mixed-monomial}
\begin{aligned}
& v_jv_k \in I + F_4\C[[V]] \text{ for $j \neq k$}, \\
& v_j^6 \in I \cdot F_2\C[[V]] + F_8\C[[V]].
\end{aligned}
\end{equation}
Start with $W_7 = \alpha^0$. Because of its symmetry properties, this can contain no pure monomials $v_k^7$ or $v_k^8$. From the first part of \eqref{eq:mixed-monomial} one sees that $W - W_7 \in I \cdot F_5\C[[V]] + F_9\C[[V]]$. By appropriately choosing $f_{5,1},\dots,f_{5,3} \in F_5\C[[V]]$, one can achieve that
\begin{equation}
\begin{aligned}
& W_7(v_1+f_{5,1}(v),\dots,v_3+f_{5,3}(v)) \\ & \equiv W_7 + f_{5,1} \partial_1 W + \cdots + f_{5,3} \partial_3 W \;\;\mathrm{mod}\; F_{9}\C[[V]]\\
&\equiv W \;\;\mathrm{mod}\; F_{9}\C[[V]].
\end{aligned}
\end{equation}
Here, the error term includes the differences $f_{5,k}(\partial_k W_7 - \partial_k W)$, as well as quadratic and higher terms in the Taylor expansion, all of which lie in $F_{11}\C[[V]]$. Moreover, one can a posteriori average the coordinate change to make it suitably equivariant. The result is a function $W_9 \equiv W$ mod $F_9\C[[V]]$, with the same symmetry properties as $W_7$ itself.

From then on, one uses a slight variant of the same strategy. Suppose that for some odd $r \geq 9$ we have a function $W_r \equiv W$ mod $F_r\C[[V]]$, which is odd and $G$-invariant. By \eqref{eq:mixed-monomial} one can write $W - W_r \in I \cdot F_{r-4}\C[[V]] + F_{r+2}\C[[V]]$. By appropriately choosing $f_{r-4,1},\dots,f_{r-4,3} \in F_{r-4}\C[[V]]$, one can achieve that $W_r(v_1 + f_{r-4,1}(v),\dots,v_3 + f_{r-4,3}(v)) \equiv W$ mod $F_{r+2}\C[[V]]$. After averaging this coordinate transformation to make it equivariant, one gets a function $W_{r+2}$ which can be used in the next step. This process yields an infinite sequence of coordinate changes, which are of increasingly high order, hence whose infinite composition converges. Alternatively, one can break off after a few steps and apply an equivariant version of the general finite determinacy theorem (see the references given above).

The conclusion is that, after acting by the exponential of some element $\gamma = (\gamma^1,0) \in \g^0$, we may assume that our Maurer-Cartan element is of the form $(W,\alpha^2)$, where $\alpha^2 \in F_4\C[[V]] \otimes \Lambda^2(V)$. Finite-dimensionality of $\C[[V]]/I$ implies that the $\partial_kW$ form a regular sequence in $\C[[V]]$, which in turn implies that the complex \eqref{eq:koszul-complex} is a resolution of $\C[[V]]/I$ \cite[Corollary 4.5.5]{weibel}. Hence $\alpha^2 = -\iota_{dW}\gamma^3$, and again one can choose $\gamma^3$ to be odd and $G$-invariant. Moreover, by looking at the low degrees in the Taylor expansion, it follows that $\gamma^3 \in F_2\C[[V]] \otimes \Lambda^3(V)$, hence lies in $\g$. According to \eqref{eq:gamma3-action}, the action of the exponential of $(0,\gamma^3)$ transforms $(W,\alpha^2)$ into $(W,0)$, which completes the proof of Lemma \ref{th:singularity}.

\section{A classification theorem\label{sec:classify}}

We correspondingly modify the noncommutative geometry side. Take $V = \C^3$ and $A = \Lambda(V)$, with the same $G \subset \mathit{SL}(V) = \mathrm{Aut}(A)$ as before. Define a graded vector space $\h$ by
\begin{equation}
\h^d \;\; = \!\!\!\!\! \prod_{\substack{3i+j-4k = 3d+3 \\ k \geq 0, \, i \geq d+2}}\!\!\!\!\! \Hom^j(A^{\otimes i},A)^G \, \hbar^k. \label{eq:h-algebra}
\end{equation}
The parity of the grading agrees with the one previously used in our discussion for $\mathit{CC}(A,A)$. Hence, the Hochschild differential and Gerstenhaber bracket turn $\h$ into a dg Lie algebra. It is filtered nilpotent, with $L_r\h^d$ being the part of \eqref{eq:h-algebra} where $i \geq d+1+r$.

\begin{lemma} \label{th:kontsevich-2}
There is a filtered $L_\infty$-quasi-isomorphism $\Phi: \h \rightarrow \g$ whose first term $\Phi^1$ is (the obvious $\hbar$-linear extension of) the Hochschild-Kostant-Rosenberg map \eqref{eq:hkr}.
\end{lemma}

This is a direct consequence of Theorem \ref{th:kontsevich}. It is useful to think of $\h$ as obtained from $\mathit{CC}(A,A)$ by the following process. One starts with $\mathit{CC}(A,A)[[\hbar]]$ with its traditional grading, and then modifies that to a $\Q$-grading by giving $\hbar$ degree $4/3$, and subtracting $2/3$ times the weight of the action of the central $\C^* \subset \mathit{GL}(V)$. Restrict to the subspace where this $\Q$-grading is integral, and where the weight of the $\C^*$-action is strictly less than four times the order of $\hbar$. In terms of \eqref{eq:h-algebra} the latter condition says that $j < 4k$, which is equivalent to $i > d+1$. Finally, take the $G$-invariant part. It is easy to check that the same process produces $\g$ from $\C[[V]] \otimes \Lambda(V)$. Since the terms $\Phi^r: \mathit{CC}(A,A)[[\hbar]]^{\otimes r} \rightarrow \C[[V]] \otimes \Lambda(V)[[\hbar]]$ of the Kontsevich $L_\infty$-homomorphism are $\mathit{GL}(V)$-equivariant and respect powers of $\hbar$, they restrict to maps $\h^{\otimes r} \rightarrow \g$ of the correct degree $1-r$. The filtrations can be defined in similar terms, showing that $\Phi$ satisfies \eqref{eq:phi-filtered}. Moreover, since parities of gradings are preserved, the restrictions satisfy the necessary symmetry and $L_\infty$-homomorphism conditions.

We will use Lemma \ref{th:kontsevich-2} to transfer classification problems for Maurer-Cartan solutions from $\h$ to $\g$. However, before bringing this theory to bear, let's look at the meaning of such solutions. A general $\alpha \in \h^1$ consists of $i$-linear components $\alpha^i$ for $i \geq 3$, each of which is in turn of the form $\alpha^i = \alpha^i_0 + \hbar \alpha^i_1 + \cdots$, with
\begin{equation}
\alpha^i_k \in \Hom^{6-3i+4k}(A^{\otimes i},A)^G.
\end{equation}
Note that for each fixed $i$, there are only finitely many $k$ such that $\alpha^i_k \neq 0$, for degree reasons. Define multilinear maps $\mu^i: A^{\otimes i} \rightarrow A$ whose $\Z/2$-grading is $i$ by setting
\begin{equation} \label{eq:wedge}
\begin{cases}
 \mu^1 = 0, \\
 \mu^2(a_2,a_1) = (-1)^{|a_1|} a_2 \wedge a_1, \\
 \mu^i = \alpha^i_0 + \alpha^i_1 + \cdots \text{ for $i \geq 3$.}
\end{cases}
\end{equation}
In parallel with our previous general discussion, $\alpha$ is a solution of the Maurer-Cartan equation iff $\mu$ is a $\Z/2$-graded $A_\infty$-structure on $A$. Of course, this structure is automatically $G$-invariant as well. Next, suppose that we have two solutions of the Maurer-Cartan equation, related by the exponentiated action of some $\gamma \in \h^0$. Then, the the associated $A_\infty$-structures $\mu$, $\tilde\mu$ are related by a $G$-equivariant $\Z/2$-graded $A_\infty$-isomorphism $\phi$, whose first term is $\phi^1 = \id$. One gets $\phi$ from $\gamma$ by the formulae from \eqref{eq:exp-gamma}, with the simplifications coming from $\gamma^1 = 0$.

Let's look explicitly at some of the simplest terms which $\alpha$ and $\gamma$ can have. First of all,
\begin{equation} \label{eq:g-zero}
\begin{aligned}
& \Hom^1(A^{\otimes 3},A)^G = 0, \\
& \Hom^{-2}(A^{\otimes 3},A)^G = 0.
\end{aligned}
\end{equation}
The first part of this implies that $\alpha^3_1 = 0$. In view of that, the simplest nontrivial components of the Maurer-Cartan equation are
\begin{equation} \label{eq:low-order}
\begin{aligned}
& \partial \alpha^3_0 = 0, \\
& \partial \alpha^4_1 = 0, \\
& \partial \alpha^5_1 + [\alpha_0^3,\alpha_1^4] = 0.
\end{aligned}
\end{equation}
$\alpha^3_0$ is a cocycle, which under the Hochschild-Kostant-Rosenberg map goes to
\begin{equation} \label{eq:alpha-0-3-image}
\Phi^1(\alpha^3_0) \in \Sym^3(V^\vee)^G = \C \cdot v_1 v_2 v_3.
\end{equation}
The next term $\alpha_4^1$ is again a cocycle, whose cohomology class is determined by
\begin{equation} \label{eq:alpha-1-4-image}
\begin{aligned}
\Phi^1(\alpha_1^4) & \in \big(\Sym^4(V^\vee) \otimes \Lambda^2(V)\big)^G \\
& = \C \cdot (v_1^4 \otimes \xi_2 \wedge \xi_3) \oplus \C \cdot (v_2^4 \otimes \xi_3 \wedge \xi_1)
\oplus \C \cdot (v_3^4 \otimes \xi_1 \wedge \xi_2).
\end{aligned}
\end{equation}
Since $\Phi^1$ induces a Lie algebra homomorphism in cohomology, it follows from \eqref{eq:low-order} that the Nijenhuis bracket $[\Phi^1(\alpha^3_0),\Phi^1(\alpha_4^1)] = 0$ vanishes, but that is only possible if one of the two classes involved is zero. The case of interest to us, which we will concentrate on from now on, is when $\Phi^1(\alpha^3_0) \neq 0$, which means that $\Phi^1(\alpha_4^1) = 0$. From the second part of \eqref{eq:g-zero}, it follows a fortiori that there are no Hochschild coboundaries in $\Hom^{-2}(A^{\otimes 4},A)^G$. Hence, vanishing of the cohomology class of $\alpha_4^1$ means that the cocycle itself is zero. That in turns means that $\alpha^5_1$ is itself a cocycle, with cohomology class
\begin{equation} \label{eq:alpha-1-5-image}
\Phi^1(\alpha^5_1) \in \Sym^5(V^\vee)^G = \C \cdot v_1^5 \oplus \C \cdot v_2^5 \oplus \C \cdot v_3^5.
\end{equation}

Next, let's analyze the action of $\gamma \in \h^0$. Both $\gamma_1^2 \in \Hom^1(A^{\otimes 2},A)^G$ and $\gamma_1^3 \in \Hom^{-2}(A^{\otimes 3},A)^G$ vanish by \eqref{eq:g-zero}. Assuming as before that $\Phi^1(\alpha^3_0) \neq 0$, the infinitesimal gauge transformation $\beta = - \partial\gamma + [\gamma,\alpha]$ has components
\begin{equation}
\begin{aligned}
 & \beta_0^3 = -\partial \gamma_0^2, \\
 & \beta_1^4 = 0, \\
 & \beta_1^5 = -\partial \gamma_1^4.
\end{aligned}
\end{equation}
Hence, the cohomology classes \eqref{eq:alpha-0-3-image} and \eqref{eq:alpha-1-5-image} are preserved, which means that they are invariants of the equivalence class of the Maurer-Cartan element.

\begin{proposition} \label{th:classification}
Up to equivalence, there is a unique Maurer-Cartan element $\alpha \in \h^1$ such that
$\Phi^1(\alpha^3_0) = - v_1v_2v_3$ and $\Phi^1(\alpha^5_1) = v_1^5 + v_2^5 + v_3^5$.
\end{proposition}

By construction, $\Phi$ is filtered and preserves powers of $\hbar$. Recall from the previous computation that $\alpha^3_1 = \alpha^4_1 = 0$. Moreover, $\Phi^r(\alpha^{i_r}_0,\dots,\alpha^{i_1}_0)$ vanishes for degree reasons (thinking back to the original grading as in Theorem \ref{th:kontsevich}), with the single exception of $\Phi^1(\alpha^3_0)$. Hence, the leading terms in the image of $\alpha$ under \eqref{eq:map-alpha} are
\begin{equation}
\begin{aligned}
\tilde\alpha & = \Phi^1(\alpha) + \half \Phi^2(\alpha,\alpha) + \cdots \\
& \equiv \Phi^1(\alpha^3 + \alpha^4 + \alpha^5) + \half \Phi^2(\alpha^3,\alpha^3) + \Phi^2(\alpha^3,\alpha^4)
\; \text{mod} \; L_4\g^1 \\
& \equiv \Phi^1(\alpha^3_0) + \hbar \Phi^1(\alpha^5_1) \; \text{mod} \; L_4\g^1 + (\hbar^2 \g)^1.
\end{aligned}
\end{equation}
Closer inspection shows that $L_4\g^1$ and $(\hbar^2\g)^1$ coincide. Either of them consists of those elements $\tilde\alpha = (\tilde\alpha^0,\tilde\alpha^2) \in \g^1$ such that $\tilde\alpha^0 \in F_7\C[[V]]$ and $\tilde\alpha^2 \in F_6\C[[V]] \otimes \Lambda^2(V)$. This means that $\tilde\alpha^0 \equiv W$ mod $F_7\C[[V]]$. Lemma \ref{th:singularity} shows that such a Maurer-Cartan solution is unique up to equivalence. Hence,  the same holds for the original solution $\alpha$, by Lemma \ref{th:filtered}.

Even though that is basically a repetition, it may still make sense to reformulate the outcome in a way which is more directly relevant to applications. Let $\mu$ be a $G$-equivariant $\Z/2$-graded $A_\infty$-structure on $A$, where $\mu^1 = 0$, $\mu^2$ is the ordinary product up to sign changes as in \eqref{eq:wedge}, and where the higher order structures can be written as sums of components $\mu^i_k$ of degree $6-3i+4k$. Assume moreover that for $\xi \in V \subset A$,
\begin{equation} \label{eq:higher-35}
\mu^3_0(\xi,\xi,\xi) = -\xi_1\xi_2\xi_3, \quad \mu^5_1(\xi,\xi,\xi,\xi,\xi) = \xi_1^5 + \xi_2^5 + \xi_3^5.
\end{equation}
These requirements determine $\mu$ uniquely up to $G$-equivariant $A_\infty$-iso\-mor\-phisms. From now on, we will generally write $\A$ for any $A_\infty$-algebra belonging to this isomorphism class.

\section{General aspects of the Fukaya category\label{sec:fukaya}}

Temporarily, allow $M$ to be a closed connected oriented surface of any genus $\geq 2$. The Fukaya category $\F(M)$ is a $\Z/2$-graded $A_\infty$-category over $\C$. We'll begin by giving a description of this category on the cohomological level, and then discuss some properties which can be stated independently of the more technical aspects of the chain level construction.

Let $S(\mathit{TM}) \rightarrow M$ be the tangent circle bundle; this can be defined without reference to a metric, as the bundle of oriented real lines in $\mathit{TM}$. Choose a symplectic form $\omega$ on $M$, and a one-form $\theta$ on $S(\mathit{TM})$ whose exterior derivative is the pullback of $\omega$; this exists because the tangent bundle has nonzero degree. Consider connected Lagrangian submanifolds, which are of course just simple closed curves $L \subset M$. Let $\sigma: L \rightarrow S(\mathit{TM})|L$ be the section given by the tangent spaces of $L$, for some choice of orientation. We say that $L$ is balanced if $\int_L \sigma^*\theta = 0$; this property is independent of the orientation, since the sections $\pm \sigma$ are fibrewise homotopic.

\begin{remark}
A nullhomologous curve is balanced iff it divides $M$ into halves $M_\pm$ such that $\chi(M_+)/\mathrm{area}(M_+) = \chi(M_-)/\mathrm{area}(M_-)$. Thus, contractible curves can never be balanced. Every other isotopy class of curves contains a balanced representative, which is unique up to Hamiltonian isotopy.
\end{remark}

Objects of $\F(M)$ are balanced curves $L$ equipped with orientations and $\mathit{Spin}$ structures. On the cohomology level, the morphisms
\begin{equation}
H(\mathrm{hom}_{\F(M)}(L_0,L_1)) = \mathit{HF}^*(L_0,L_1)
\end{equation}
are the Lagrangian Floer cohomology groups. In particular, for any object $L$ we have a canonical isomorphism $H(\mathrm{hom}_{\F(M)}(L,L)) = \mathit{HF}^*(L,L) \iso H^*(L;\C)$. Still remaining on the cohomology level, composition of morphisms is given by Donaldson's holomorphic triangle product. In particular, the isomorphism $\mathit{HF}^*(L,L) \iso H^*(L;\C)$ is compatible with the ring structure. Moreover, any two balanced curves which are isotopic (compatibly with the orientations and {\em Spin} structures) give rise to isomorphic objects of $H^0(\F(M))$.

\begin{remark}
In the case of the torus \cite{polishchuk-zaslow98}, the Fukaya category is defined over a Novikov field, which is a field of formal Laurent series with a parameter $t$. If one wants to define a Fukaya category of a higher genus surface containing {\em all} Hamiltonian isotopy classes of curves as objects, Novikov fields appear there as well. However, if one then restricts attention to balanced curves, all resulting series turn out to be finite (Laurent polynomials), hence one can set $t = 1$ and work over $\C$, which is what we are doing here. This phenomenon, usually called monotonicity, is familiar to symplectic geometers (it appears in the literature mainly in the context of Floer cohomology for Fano manifolds; see \cite{oh93}, and for a formulation closer to the one adopted here, \cite[Remark 3.1.4]{wehrheim-woodward06}).
\end{remark}

The definition of the objects in the Fukaya category involves $\theta$, or rather its equivalence class modulo exact one-forms (compare \cite[Appendix A]{abouzaid08}). This has some consequences for functoriality. Namely, suppose that we have two choices of symplectic forms $\omega$, $\tilde\omega$ and correspondingly one-forms $\theta$, $\tilde\theta$. Given a symplectomorphism $\phi: (M,\omega) \rightarrow (M,\tilde\omega)$, we can consider the induced map $S(D\phi): S(\mathit{TM}) \rightarrow S(\mathit{TM})$, which defines a class
\begin{equation}
[S(D\phi)^*\tilde\theta - \theta] \in H^1(S(\mathit{TM});\R).
\end{equation}
We say that $\phi$ is balanced if this class vanishes, in which case it induces a quasi-isomorphism between the associated Fukaya categories. Every connected component of the space of symplectomorphisms contains such representatives, and they are unique up to Hamiltonian isotopy. This, together with invariance under rescaling of $\omega$ and $\theta$, implies that $\F(M)$ is independent of the additional choices up to quasi-isomorphisms. Hence, it is justifiable to talk of ``the Fukaya category of $M$''. For the same reason, the mapping class group of $M$ acts on $\F(M)$.

Before continuing, we need to recall a few homological algebra notions from \cite{keller99} or \cite[Sections 3--4]{seidel04}. Given any $\Z/2$-graded $A_\infty$-category $\B$, one can consider the associated dg category $\mathit{mod}(\B)$ of right $A_\infty$-modules. This comes with a canonical cohomologically full and faithful $A_\infty$-functor $\B \rightarrow \mathit{mod}(\B)$, the Yoneda embedding. The associated cohomology level category $\mathit{Mod}(\B) = H^0(\mathit{mod}(\B))$ is a triangulated category, with the property that the twofold shift functor is isomorphic to the identity, and moreover it is split-closed (also called Karoubi complete; this means that any idempotent endomorphism of an object leads to a splitting of that object as a direct sum). Take the smallest full subcategory of $\mathit{Mod}(\B)$ which contains the image of the Yoneda embedding, is triangulated, and split-closed. We denote this by $D^\pi(\B)$, and call it the split-closed derived category of $\B$. %Note that by definition, this has an underlying dg category.

\begin{remark}
The only point in the above discussion which might not be entirely familiar is the fact that $\mathit{Mod}(\B)$ is split-closed. A short proof goes as follows. The dg category $\mathit{mod}(\B)$ itself comes with a Yoneda embedding, which induces a functor $\mathit{Mod}(\B) \rightarrow \mathit{Mod}(mod(\B))$. In the other direction, we have a restriction functor $\mathit{Mod}(\mathit{mod}(\B)) \rightarrow \mathit{Mod}(\B)$; and the composition of the two is the identity on $\mathit{Mod}(\B)$. From \cite[Section 4b]{seidel04} we know that any idempotent endomorphism in $\mathit{Mod}(\B)$ leads to a splitting of the associated Yoneda module in $\mathit{Mod}(\mathit{mod}(\B))$. The image of that under restriction is the desired splitting in $\mathit{Mod}(\B)$.
\end{remark}

Let $\A \subset \B$ be a full $A_\infty$-subcategory. In that case, $D^\pi(\A)$ is canonically equivalent to the smallest split-closed triangulated full subcategory of $D^\pi(\B)$ containing all objects of $\A$. Objects of $\B$ which, up to isomorphism, lie in $D^\pi(\A)$, are said to be split-generated by the objects of $\A$. If this holds for all of $\B$, which means that the embedding $D^\pi(\A) \rightarrow D^\pi(\B)$ is an equivalence, we say that the objects of $\A$ split-generate $\B$.

To apply this to the Fukaya category, we need to recall some facts about the action of Dehn twists. Let $L_0,L_1$ be objects of $\F(M)$, where the $\mathit{Spin}$ structure on $L_1$ is nontrivial. The Dehn twist $\tau_{L_1}$ is a balanced symplectic automorphism of $M$, hence $\tau_{L_1}(L_0)$ is again a balanced curve. We then have an exact triangle in $D^\pi(\F(M))$ of the form
\begin{equation} \label{eq:exactness}
\xymatrix{
\mathit{HF}^*(L_1,L_0) \otimes L_1 \ar[r]^-{\mathrm{ev}} & L_0 \ar[d] \\
& \tau_{L_1}(L_0). \ar[ul]^{[1]}
}
\end{equation}
Here, $\mathit{HF}^*(L_1,L_0) \otimes L_1$ is a direct sum of copies of $L_1$ and its shifted version $L_1[1]$, with one summand for each generator of the Floer cohomology group, and $\mathrm{ev}$ is the canonical evaluation map. An equivalent statement is that $\tau_{L_1}(L_0)$ is isomorphic to the cone of $\mathrm{ev}$. The construction of the exact triangle runs parallel to \cite[Corollary 17.18]{seidel04}, which means that it is based on a version of the long exact sequence from \cite{seidel01}.

\begin{lemma} \label{th:karoubi-1}
Let $\{L_1,\dots,L_r\}$, $r \geq 1$, be objects of $\F(M)$ whose $\mathit{Spin}$ structures are nontrivial. Let $L_0$ be another object, such that $\tau_{L_r}\cdots\tau_{L_1}(L_0)$ is isotopic to $L_0$ with the orientation reversed. Then $L_0$ is split-generated by $\{L_1,\dots,L_r\}$.
\end{lemma}

This is similar to \cite[Proposition 19.7]{seidel04}. The composition of the vertical arrows in \eqref{eq:exactness} yields a morphism
\begin{equation} \label{eq:ttt}
L_0 \longrightarrow \tau_{L_r}\cdots\tau_{L_1}(L_0) \iso L_0[1].
\end{equation}
This morphism is an element of $\mathit{HF}^1(L_0,L_0) \iso H^1(L_0;\C)$, hence its square is automatically zero. This is precisely what's needed to make the argument in \cite[p.\ 70]{seidel04} go through.

\begin{lemma} \label{th:karoubi-2}
Let $\{L_1,\dots,L_r\}$ be objects of $\F(M)$ whose $\mathit{Spin}$ structures are nontrivial, and such that  $\tau_{L_r}\cdots\tau_{L_1}$ is isotopic to the identity. Then they split-generate $\F(M)$.
\end{lemma}

The basic strategy is the same, but with an additional geometric step. Take an arbitrary $L_0$, and consider the analogue of \eqref{eq:ttt}, which this time is an element of $\mathit{HF}^0(L_0,L_0) \iso H^0(L_0;\C)$. By construction of the exact triangle, this element admits the following description. Consider the Lefschetz fibration with fibre $M$ and vanishing cycles $\{L_1,\dots,L_r\}$. Fix a generic almost complex structure which makes the fibration map pseudo-holomorphic, and consider the associated moduli space of pseudo-holomorphic sections. This may have components of different dimensions, but (due to the balancing condition, and the fact that the fibres contain no holomorphic spheres) the component of any fixed dimension is compact. By considering the evaluation map at a point, as in Gromov-Witten theory, one gets an even-dimensional cohomology class in $M$, which we call the section class. For any $L_0$, the morphism \eqref{eq:ttt} is the image of the section class under the restriction map
\begin{equation} \label{eq:restrict-sigma}
H^0(M;\C) \rightarrow H^0(L_0;\C).
\end{equation}
Assume that the section class has a nontrivial component in $H^0(M;\C)$. This means that through every point of $M$ there is a pseudo-holomorphic section with zero selfintersection. Standard methods from four-dimensional symplectic topology \cite{mcduff90b} then ensure that these sections foliate the total space of our Lefschetz fibration, which is a contradiction. Hence the image of the section class under \eqref{eq:restrict-sigma} vanishes, allowing one to proceed as before.

\begin{remark}
The section class itself is not zero in general. For instance, if the fibration is constructed by blowing up a Lefschetz pencil, every base point of the pencil gives rise to a section, which contributes $1$ to the $H^2(M;\C)$ component of the section class (however, taking the fibre connect sum of the fibration with itself corresponds to passing to the cup-square of the section class, which will kill it).
\end{remark}

\section{Technical aspects of the Fukaya category\label{sec:technical}}

In this section, we take a closer look at the definition of the Fukaya category. Since our target space is a surface, the naive idea is that after appealing to the uniformization theorem, the $A_\infty$-structure maps should be computable purely combinatorially by counting polygons. This is true in many cases but fails to hold in general, due to transversality issues, which any proper definition must address. There are several approaches, all of which are ultimately equivalent (meaning that they give rise to different but quasi-equivalent $A_\infty$-categories). We follow the ``Morse-Bott'' type approach, in a version which borrows some aspects of \cite{schwarz96} and \cite{cornea-lalonde06}.

Fix a countable set $\scrL$ of balanced curves on $M$ with the following properties. Each nontrivial isotopy class has at least one representative in $\scrL$. Moreover, any two distinct curves in $\scrL$ intersect transversally, and any three distinct curves have no common point. From now on, when defining $\F(M)$, we will only allow curves taken from $\scrL$ (this is a technical contrivance, which is ultimately irrelevant: any two choices of $\scrL$ lead to quasi-equivalent $A_\infty$-categories).

Suppose that $L_0$ and $L_1$ are objects, and that the underlying curves are distinct, hence transverse. In this case, the morphism space between them is the unperturbed Floer cochain complex
\begin{equation} \label{eq:cf}
\mathrm{hom}_{\F(M)}(L_0,L_1) = \mathit{CF}^*(L_0,L_1) = \bigoplus_x \C x,
\end{equation}
where the sum is over all intersection points $x \in L_0 \cap L_1$. The generator associated to $x$ is even if the local intersection number is $-1$, and odd otherwise. Next, suppose that $(L_0,\dots,L_d)$ is a collection of objects, whose underlying curves are pairwise different. In this case, the coefficients of the $A_\infty$-composition
\begin{equation}
\begin{aligned}
& \mu^d: \mathit{CF}^*(L_{d-1},L_d) \otimes \cdots \otimes \mathit{CF}^*(L_0,L_1) \longrightarrow \mathit{CF}^*(L_0,L_d), \\
& \mu^d(x_d,\dots,x_1) = \sum_{x_0} m(x_0,\dots,x_d) x_0
\end{aligned}
\end{equation}
are numbers $m(x_0,\dots,x_d) \in \Z$ obtained by a signed count of immersed polygons. The construction is well-known. We describe it briefly, and refer to \cite[Section 13]{seidel04} for details; for other versions see \cite{polishchuk98,chekanov99,auroux-katzarkov-orlov04}.

Fix a Riemann surface structure on $M$, compatible with its symplectic orientation. Let $(\zeta_0,\dots,\zeta_d)$, $d \geq 1$, be distinct boundary points on the closed disc $D \subset \C$, ordered in accordance with the boundary orientation. Consider holomorphic maps $u: D \setminus \{\zeta_0,\dots,\zeta_d\} \rightarrow M$ which map the boundary sides to $(L_0,\dots,L_d)$, and which extend continuously to $D$, taking each $\zeta_k$ to $x_k$. Each $(\zeta_0,\dots,\zeta_d,u)$ has a virtual dimension, and we consider only ones of virtual dimension zero. The equivalence relation is that $(\zeta_0,\dots,\zeta_d,u) \sim (\tilde\zeta_0,\dots,\tilde\zeta_d,\tilde{u})$ if there is an automorphism $\phi: D \rightarrow D$ such that $\phi(\zeta_k) = \tilde{\zeta}_k$ and $u = \tilde{u} \circ \phi$. Denote the resulting space of equivalence classes by $\scrM(x_0,\dots,x_d)$. Since constant maps $u$ are excluded by our assumptions, automatic regularity \cite[Lemma 13.2]{seidel04} ensures that each point of this space is regular, hence contributes $\pm 1$ to $m(x_0,\dots,x_d)$. To translate this into combinatorics, one observes first that points of $\scrM(x_0,\dots,x_d)$ correspond bijectively to immersed polygons with sides on the $L_k$, and corners at $x_k$. To compute the sign with which each polygon contributes, we pick, for each $L_k$ such that the $\mathit{Spin}$ structure is nontrivial, a marked point $\circ_k \in L_k$ which is not an intersection point with any of the other curves in $\scrL$, as well as a trivialization of the $\mathit{Spin}$ structure away from that point. If the $L_k$ are oriented as in Figure \ref{fig:or}(i), and if none of the points $\circ_k$ lie on the boundary of our immersed polygon, its contribution to $m(x_0,\dots,x_d)$ is $+1$. The general rule is obtained from this by the following sign changes. The orientation of $L_0$ is irrelevant. Reversing the orientation of $L_k$, $0<k<d$, changes the sign by $(-1)^{|x_k|}$. Reversing the orientation of $L_d$ changes the sign by $(-1)^{|x_0|+|x_d|}$. Finally, for every time that the boundary of the polygon passes over one of the points $\circ_k$, for $0 \leq k \leq d$, we change the sign by $(-1)$.
\begin{figure}[ht]
\begin{centering}
\begin{picture}(0,0)%
\includegraphics{or.pstex}%
\end{picture}%
\setlength{\unitlength}{3552sp}%
\begingroup\makeatletter\ifx\SetFigFont\undefined%
\gdef\SetFigFont#1#2#3#4#5{%
  \reset@font\fontsize{#1}{#2pt}%
  \fontfamily{#3}\fontseries{#4}\fontshape{#5}%
  \selectfont}%
\fi\endgroup%
\begin{picture}(4977,2274)(136,-4123)
\put(1276,-3661){\makebox(0,0)[lb]{\smash{{\SetFigFont{11}{13.2}{\rmdefault}{\mddefault}{\updefault}{\color[rgb]{0,0,0}$L_1$}%
}}}}
\put(151,-2236){\makebox(0,0)[lb]{\smash{{\SetFigFont{11}{13.2}{\rmdefault}{\mddefault}{\updefault}{\color[rgb]{0,0,0}(i)}%
}}}}
\put(2851,-2236){\makebox(0,0)[lb]{\smash{{\SetFigFont{11}{13.2}{\rmdefault}{\mddefault}{\updefault}{\color[rgb]{0,0,0}(ii)}%
}}}}
\put(4201,-3286){\makebox(0,0)[lb]{\smash{{\SetFigFont{11}{13.2}{\rmdefault}{\mddefault}{\updefault}{\color[rgb]{0,0,0}$L_1$}%
}}}}
\put(4351,-2386){\makebox(0,0)[lb]{\smash{{\SetFigFont{11}{13.2}{\rmdefault}{\mddefault}{\updefault}{\color[rgb]{0,0,0}$L_2$}%
}}}}
\put(2950,-2711){\makebox(0,0)[lb]{\smash{{\SetFigFont{11}{13.2}{\rmdefault}{\mddefault}{\updefault}{\color[rgb]{0,0,0}$L_0=L_3$}%
}}}}
\put(301,-3361){\makebox(0,0)[lb]{\smash{{\SetFigFont{11}{13.2}{\rmdefault}{\mddefault}{\updefault}{\color[rgb]{0,0,0}$L_0$}%
}}}}
\put(1726,-2686){\makebox(0,0)[lb]{\smash{{\SetFigFont{11}{13.2}{\rmdefault}{\mddefault}{\updefault}{\color[rgb]{0,0,0}$L_2$}%
}}}}
\put(751,-2461){\makebox(0,0)[lb]{\smash{{\SetFigFont{11}{13.2}{\rmdefault}{\mddefault}{\updefault}{\color[rgb]{0,0,0}$L_3$}%
}}}}
\put(3450,-3100){\makebox(0,0)[lb]{\smash{{\SetFigFont{11}{13.2}{\rmdefault}{\mddefault}{\updefault}{\color[rgb]{0,0,0}$e$}%
}}}}
\end{picture}%
\caption{\label{fig:or}}
\end{centering}
\end{figure}

\begin{lemma} \label{th:monotonicity-1}
Let $S$ be a compact oriented surface with boundary, and $w: S \rightarrow M$ a map which takes each boundary component of $S$ to a balanced curve. Then the normalized area $[w^*\omega]/\mathrm{area}(M)$ and the normalized relative Chern class $c_1^{\mathrm{rel}}(w^*\mathit{TM})/\chi(M)$, both lying in $H^2(S,\partial S;\R)$, agree.
\end{lemma}

This is a straightforward fact, which is worth while mentioning because it leads to the basic compactness result for the spaces $\scrM(x_0,\dots,x_d)$. This is a ``monotonicity'' style consideration, which we summarize briefly. Consider two maps $u,\tilde{u}$ which contribute to $m(x_0,\dots,x_d)$. By gluing together their domains topologically, we get a map $w: S \rightarrow M$, where $S$ is a genus zero surface with $d+1$ boundary circles, and where the images of the boundary components lie on the balanced curves $L_k$; this is unique up to homotopy within the class of such maps. An index theory argument shows that the relative Chern class of $w$ is zero, hence by Lemma \ref{th:monotonicity-1} that $\int w^*\omega$ vanishes. But by construction, that implies that the areas $\int u^*\omega$ and $\int \tilde{u}^*\omega$ are the same. From this, a compactness argument shows that $\scrM(x_0,\dots,x_d)$ is a finite set. In fact, one can translate this argument into combinatorics, where it becomes elementary.

Now consider two objects $L_0,L_1$ such that the underlying curves agree. In that case, we fix a metric and a Morse function $f_{01}$ on that curve, with a unique minimum and maximum which are both distinct from the intersection points with any other curve in $\scrL$. Denote the minimum by $e$ and the maximum by $q$. We then define the morphism space to be the Morse cochain space
\begin{equation} \label{eq:hom01}
\mathrm{hom}_{\F(M)}(L_0,L_1) = \mathit{CM}^*(f_{01}) = \C e \oplus \C q.
\end{equation}
Suppose that $L_0,L_1$ actually have the same orientation and isomorphic {\em Spin} structures. Then, the $\Z/2$-grading on \eqref{eq:hom01} coincides with the ordinary Morse index, and the differential $\mu^1$ is the Morse differential of $f_{01}$, hence vanishes. This has to be suitably modified for the other cases. For instance, suppose that the orientations agree, but that the {\em Spin} structure on $L_1$ differs from that on $L_0$ by twisting with a nontrivial double cover $\xi \rightarrow L_0$. Then, $\mu^1$ is the Morse differential with twisted coefficients in $\xi \otimes_{\Z/2} \C$, hence acyclic. See \cite[Example 13.5]{seidel04} for further discussion.

Here is a simple class of higher order compositions involving \eqref{eq:hom01}. Take objects $(L_0,\dots,L_d)$, $d \geq 2$, such that $L_{i-1} = L_i$ agree as curves for a single $0<i\leq d$, and where the underlying curves are otherwise distinct (which means that there are $d$ distinct curves among them). Denote by $f_{i-1,i}$ the function used to define $\mathrm{hom}_{\F(M)}(L_{i-1},L_i)$. Choose intersection points $x_0 \in L_0 \cap L_d$, $x_k \in L_{k-1} \cap L_k$ (for $1 \leq k \leq d$ with $k \neq i$), and a critical point $x_i$ of $f_{i-1,i}$. Consider again holomorphic maps $u: D \setminus \{\zeta_0,\dots,\zeta_d\} \rightarrow M$, but where the extension at $\zeta_i$ is now smooth, and satisfies
\begin{equation} \label{eq:unstable}
u(\zeta_i) \in W^u(x_i) \subset L_i,
\end{equation}
$W^u(x_i)$ being the unstable manifold of $x_i$ for the gradient flow. The moduli space of such maps of virtual dimension zero is again a finite set $\scrM(x_0,\dots,x_d)$, and an appropriate signed count of points in it yields integers $m(x_0,\dots,x_d) \in \Z$ which are the coefficients of the composition map
\begin{equation}
\mu^d: \mathit{CF}^*(L_{d-1},L_d) \otimes \cdots \otimes \mathit{CM}^*(f_{i-1,i}) \otimes \cdots \otimes \mathit{CF}^*(L_0,L_1) \longrightarrow \mathit{CF}^*(L_0,L_d).
\end{equation}
Again, this can be translated into combinatorics, as follows. Suppose first that $x_i = e$, where \eqref{eq:unstable} reduces to the open condition $u(\zeta_i) \neq q$. Then, the only case where $(\zeta^0,\dots,\zeta^d,u)$ has virtual dimension zero is when $d = 2$, and $u$ is the constant map at a point of $L_0 \cap L_2$. Next, consider the case $x_i = q$, where \eqref{eq:unstable} says that $u(\zeta_i) = q$. There, $m(x_0,\dots,x_d)$ can be computed by a signed count of immersed $d$-gons with an additional marked point on the appropriate boundary side, whose image is $q$.

There is another case which can be treated in the same way. Take objects $(L_0,\dots,L_d)$, $d \geq 2$, such that the first $d$ have pairwise distinct underlying curves, but that $L_0 = L_d$ as curves. Let $f_{0,d}$ be the function used to define $\mathrm{hom}_{\F(M)}(L_0,L_d)$. Choose a critical point $x_0$ of that functions, and intersection points $x_k \in L_{k-1} \cap L_k$, $k > 0$. In this case, the condition analogous to \eqref{eq:unstable} involves the stable manifold $W^s(x_0)$:
\begin{equation} \label{eq:stable}
u(\zeta_0) \in W^s(x_0) \subset L_0.
\end{equation}
Again, we have an appropriate moduli space $\scrM(x_0,\dots,x_d)$ and a signed count $m(x_0,\dots,x_d) \in \Z$, which defines the map
\begin{equation}
\mu^d: \mathit{CF}^*(L_{d-1},L_d) \otimes \cdots \otimes \mathit{CF}^*(L_0,L_1) \longrightarrow \mathit{CM}^*(f_{0,d}).
\end{equation}
Translation into combinatorics now works as follows: if $x_0 = q$, only constant triangles contribute, while for $x_0 = e$, we are counting immersed $d$-gons with an additional marked boundary point going through $e$. Having reduced the computations to combinatorics in principle, it remains to describe the signs of each polygon. We will concentrate on the cases that actually occur in our application, and give the resulting formulae without proof (verification is tedious but not difficult, following the argument from \cite{seidel04}).

{\it Constant maps:} Take two curves $L_0 \neq L_1$. The constant triangle at any point $x \in L_0 \cap L_1$ contributes to the products
\begin{equation} \label{eq:e1}
\begin{aligned}
& \mu^2(x,e),\; \mu^2(e,x): \mathit{CF}^*(L_0,L_1) \longrightarrow \mathit{CF}^*(L_0,L_1), \\
& \mu^2(x,x): \mathit{CF}^*(L_1,L_0) \otimes \mathit{CF}^*(L_0,L_1) \longrightarrow \mathit{CM}^*(f),
\end{aligned}
\end{equation}
where $f$ is the function associated to the pair $(L_0,L_0)$. In all three cases, there are no contributions from non-constant triangles (note that in the last case, we know a priori for degree reasons that the product must be a multiple of $q$). Taking signs into account, the consequence is that
\begin{equation} \label{eq:e2}
\begin{aligned}
& \mu^2(x,e) = x, \; \mu^2(e,x) = (-1)^{|x|} x, \\ & \mu^2(x,x) = (-1)^{|x|} q,
\end{aligned}
\end{equation}
where $|x| \in \Z/2$ is the degree of the generator $x \in \mathit{CF}^*(L_0,L_1)$.

{\it Non-constant maps:} Take the case $L_0 = L_d$ (where the two are assumed to have the same orientation and {\em Spin} structure). Choose intersection points $x_1,\dots,x_d$, and consider the contribution of a polygon to $m(e,x_1,\dots,x_d)$. If the orientations are as in  Figure \ref{fig:or}(ii) and our polygon avoids all the points $\circ_i$ which describe the {\em Spin} structures, then its contribution is $+1$. The rules for sign changes are the same as for Figure \ref{fig:or}(i), taking into account the fact that $x_0 = e$ is even.
\begin{figure}[hb]
\begin{centering}
\begin{picture}(0,0)%
\includegraphics{tree.pstex}%
\end{picture}%
\setlength{\unitlength}{2565sp}%
\begingroup\makeatletter\ifx\SetFigFont\undefined%
\gdef\SetFigFont#1#2#3#4#5{%
  \reset@font\fontsize{#1}{#2pt}%
  \fontfamily{#3}\fontseries{#4}\fontshape{#5}%
  \selectfont}%
\fi\endgroup%
\begin{picture}(6852,6235)(-3011,-3059)
\put(-449,2339){\makebox(0,0)[lb]{\smash{{\SetFigFont{10}{12}{\rmdefault}{\mddefault}{\updefault}{\color[rgb]{0,0,0}$L_7$}%
}}}}
\put(-1499,2264){\makebox(0,0)[lb]{\smash{{\SetFigFont{10}{12}{\rmdefault}{\mddefault}{\updefault}{\color[rgb]{0,0,0}$L_8$}%
}}}}
\put(1951,764){\makebox(0,0)[lb]{\smash{{\SetFigFont{10}{12}{\rmdefault}{\mddefault}{\updefault}{\color[rgb]{0,0,0}$L_4$}%
}}}}
\put(2326,1214){\makebox(0,0)[lb]{\smash{{\SetFigFont{10}{12}{\rmdefault}{\mddefault}{\updefault}{\color[rgb]{0,0,0}$L_5$}%
}}}}
\put(2326,1589){\makebox(0,0)[lb]{\smash{{\SetFigFont{10}{12}{\rmdefault}{\mddefault}{\updefault}{\color[rgb]{0,0,0}$L_6$}%
}}}}
\put(3151,2339){\makebox(0,0)[lb]{\smash{{\SetFigFont{10}{12}{\rmdefault}{\mddefault}{\updefault}{\color[rgb]{0,0,0}$L_7$}%
}}}}
\put(1426,1889){\makebox(0,0)[lb]{\smash{{\SetFigFont{10}{12}{\rmdefault}{\mddefault}{\updefault}{\color[rgb]{0,0,0}$L_8$}%
}}}}
\put(3826,-136){\makebox(0,0)[lb]{\smash{{\SetFigFont{10}{12}{\rmdefault}{\mddefault}{\updefault}{\color[rgb]{0,0,0}$L_3$}%
}}}}
\put(3376,-661){\makebox(0,0)[lb]{\smash{{\SetFigFont{10}{12}{\rmdefault}{\mddefault}{\updefault}{\color[rgb]{0,0,0}$L_2$}%
}}}}
\put(-749,1214){\makebox(0,0)[lb]{\smash{{\SetFigFont{10}{12}{\rmdefault}{\mddefault}{\updefault}{\color[rgb]{0,0,0}$L_5$}%
}}}}
\put(-899,1739){\makebox(0,0)[lb]{\smash{{\SetFigFont{10}{12}{\rmdefault}{\mddefault}{\updefault}{\color[rgb]{0,0,0}$L_6$}%
}}}}
\put(-2324,464){\makebox(0,0)[lb]{\smash{{\SetFigFont{10}{12}{\rmdefault}{\mddefault}{\updefault}{\color[rgb]{0,0,0}$L_1$}%
}}}}
\put(-1349,464){\makebox(0,0)[lb]{\smash{{\SetFigFont{10}{12}{\rmdefault}{\mddefault}{\updefault}{\color[rgb]{0,0,0}$L_4$}%
}}}}
\put(-824,-811){\makebox(0,0)[lb]{\smash{{\SetFigFont{10}{12}{\rmdefault}{\mddefault}{\updefault}{\color[rgb]{0,0,0}$L_2$}%
}}}}
\put(301,-211){\makebox(0,0)[lb]{\smash{{\SetFigFont{10}{12}{\rmdefault}{\mddefault}{\updefault}{\color[rgb]{0,0,0}$L_3$}%
}}}}
\put(2476,-211){\makebox(0,0)[lb]{\smash{{\SetFigFont{10}{12}{\rmdefault}{\mddefault}{\updefault}{\color[rgb]{0,0,0}$L_4$}%
}}}}
\put(3001,164){\makebox(0,0)[lb]{\smash{{\SetFigFont{10}{12}{\rmdefault}{\mddefault}{\updefault}{\color[rgb]{0,0,0}$L_4$}%
}}}}
\put(2251,2264){\makebox(0,0)[lb]{\smash{{\SetFigFont{10}{12}{\rmdefault}{\mddefault}{\updefault}{\color[rgb]{0,0,0}$L_8$}%
}}}}
\put(1501,-586){\makebox(0,0)[lb]{\smash{{\SetFigFont{10}{12}{\rmdefault}{\mddefault}{\updefault}{\color[rgb]{0,0,0}$L_1$}%
}}}}
\put(2701,-961){\makebox(0,0)[lb]{\smash{{\SetFigFont{10}{12}{\rmdefault}{\mddefault}{\updefault}{\color[rgb]{0,0,0}$L_2$}%
}}}}
\put(1276,839){\makebox(0,0)[lb]{\smash{{\SetFigFont{10}{12}{\rmdefault}{\mddefault}{\updefault}{\color[rgb]{0,0,0}$L_1$}%
}}}}
\put(3001,1664){\makebox(0,0)[lb]{\smash{{\SetFigFont{10}{12}{\rmdefault}{\mddefault}{\updefault}{\color[rgb]{0,0,0}$L_6$}%
}}}}
\put(-974,-2986){\makebox(0,0)[lb]{\smash{{\SetFigFont{10}{12}{\rmdefault}{\mddefault}{\updefault}{\color[rgb]{0,0,0}$L_1 = L_2 = L_4 = L_5$, $L_6 = L_8$}%
}}}}
\put(-974,-2761){\makebox(0,0)[lb]{\smash{{\SetFigFont{10}{12}{\rmdefault}{\mddefault}{\updefault}{\color[rgb]{0,0,0}Coincidences of the underlying curves:}%
}}}}
\put(-1274,-2086){\makebox(0,0)[lb]{\smash{{\SetFigFont{10}{12}{\rmdefault}{\mddefault}{\updefault}{\color[rgb]{0,0,0}(i)}%
}}}}
\put(2026,-2086){\makebox(0,0)[lb]{\smash{{\SetFigFont{10}{12}{\rmdefault}{\mddefault}{\updefault}{\color[rgb]{0,0,0}(ii)}%
}}}}
\end{picture}%
\caption{\label{fig:tree}}
\end{centering}
\end{figure}

Finally, we turn to the question of defining $\mu^d$, $d \geq 2$, in general, where any number of the curves $(L_0,\dots,L_d)$ may coincide. One chooses generators $(x_0,\dots,x_d)$ of the associated Floer or Morse complexes as before, but the definition of the relevant moduli space $\scrM(x_0,\dots,x_d)$ is somewhat more involved. A point of this moduli space consists of the following data (similar to, but not quite the same as, the ``clusters'' from \cite{cornea-lalonde06}):
\begin{quote}
First, we have a planar tree $T \subset \R^2$ with $d+1$ semi-infinite edges, and with vertices which are at least trivalent. The components of $\R^2 \setminus T$ should be labeled by $(L_0,\dots,L_d)$ in accordance with their natural cyclic ordering. We require that the two regions separated by any finite edge of $T$ are marked with objects which have the same underlying curve, see Figure \ref{fig:tree}(i).

Next, for every vertex $v$ of $T$ of valency $|v|$, we have a collection of marked points $\{\zeta_0^v,\dots,\zeta_{|v|-1}^v\}$ on the boundary of the disc, and a holomorphic map $u_v: D \setminus \{\zeta_0^v,\dots,\zeta_{|v|-1}^v\} \rightarrow M$, extending continuously over $D$. The boundary conditions for that map are given by the labels $L_k$ in components of $\R^2 \setminus T$ adjacent to $v$. See Figure \ref{fig:tree}(ii) for an illustration.

Now take a finite edge of $T$. Associated to this edge are its two endpoints $v_\pm$, and also marked points $\zeta_{k_\pm}^{v_\pm}$. By assumption, the two objects $(L_i,L_j)$ labeling the components of $\R^2 \setminus T$ adjacent to our edge share the same underlying curve. Let $f_{i,j}$ be the Morse function used to construct $\mathrm{hom}_{\F(M)}(L_i,L_j)$. We then ask that there should be a gradient flow line of $f_{i,j}$ of some finite nonzero length, which goes from $u(\zeta_{k_-}^{v_-})$ to $u(\zeta_{k_+}^{v_+})$. The implicit notational convention here is indicated by the direction of the arrows in Figure \ref{fig:tree}(ii).

Finally, consider an infinite edge of $T$, with its unique associated endpoint $v$ and marked point $\zeta^v_k$. If the curves labeling the two components of $\R^2 \setminus T$ adjacent to our edge are distinct, we ask that $u(\zeta_k^v)$ should lie at the relevant intersection point $x_i$. Otherwise, we impose conditions as in \eqref{eq:stable}, \eqref{eq:unstable}.
\end{quote}

It is easy to see that this generalizes the previous discussion: in all cases we had considered before, the requirements only allow the star-shaped tree $T$ (with a single vertex, hence no finite edges). Unfortunately, in general the spaces $\scrM(x_0,\dots,x_d)$ are not regular (due to the failure of gradient flow lines to intersect transversally, and to the appearance of constant holomorphic maps which have excess dimension). Hence, one has to perturb this initial definition either virtually, which leads to the construction of appropriate virtual fundamental chains on the compactifications $\bar\scrM(x_0,\dots,x_d)$ \cite{fooo}, or else by perturbing the gradient flow equations and holomorphic map equations themselves, in the manner of \cite{seidel04}. Fortunately, the only case of this more complicated formalism which we need to determine explicitly is the product structure on $\mathrm{hom}_{\F(M)}(L,L)$ for a single object $L$, which is given by
\begin{equation} \label{eq:e3}
\mu^2(e,e) = e,\; \mu^2(q,e) = q, \; \mu^2(e,q) = -q.
\end{equation}
In this specific case, this also follows from the general fact that the product reproduces the ordinary cup product on $H^*(L;\C)$.

\begin{remark}
We want to quickly mention some other definitions of the Fukaya category. \cite{fooo} also uses a Morse-Bott method, but where singular cohomology replaces Morse cohomology. In contrast, \cite{seidel04} uses Hamiltonian perturbations of the holomorphic map equation to treat $\mathrm{hom}_{\F(M)}(L_0,L_1)$ on the same footing for all pairs $(L_0,L_1)$. The equivalence of any two approaches can be proved by constructing a ``mixed'' Fukaya category which contains two copies of each object, to which the two different methods are applied; compare the discussion in \cite[Section 10a]{seidel04}. Strictly speaking, the only substantial property of the Fukaya category which we have borrowed from the literature is the existence of exact triangles \eqref{eq:exactness}, which quotes \cite{seidel04}. However, the argument leading to those triangles involves only the objects $(L_0,L_1,\tau_{L_1}(L_0))$, which moreover can be perturbed to be in general position. In that form, it carries over easily to any other framework, such as the one adopted here.
\end{remark}

\section{Gradings}

The lack of an integer grading on $\F(M)$ is unavoidable, since it is directly related to the nonvanishing of $c_1(M)$. Nevertheless, one can partially improve the situation by thinking of $c_1(M)$ as being supported at finitely many points. Namely, let $\eta^r$ be a nonzero meromorphic section of the $r$-th power of the canonical bundle $\mathit{T^*\!M}^{\otimes r}$, for $r \neq 0$. Let $D \subset M$ be the set of its zeros and poles, with the order of vanishing written as $\mathrm{ord}(\eta^r,z) \in \Z$ (a negative order signifies a pole). For every oriented $L \subset M \setminus D$ we get a map $L \rightarrow S^1$, defined by
\begin{equation} \label{eq:phase}
x \longmapsto \frac{\eta^r(X^{\otimes r})}{\|\eta^r(X^{\otimes r})\|}
\end{equation}
where $X \in \mathit{TL}_x$ is nonzero and points in positive direction. An $1/r$-grading of $L$ is a lift $a: L \rightarrow \R$ of this map. Let $\F(M,D)$ be a version of the Fukaya category, whose objects are Lagrangian submanifolds as before, with the added condition that they lie in $M \setminus D$, and moreover should come equipped with $1/r$-gradings. Nothing else changes, which in particular means that there is a full and faithful $A_\infty$-functor $\F(M,D) \rightarrow \F(M)$.

In the presence of $1/r$-gradings, the generators $x$ of $\mathrm{hom}_{\F(M)}(L_0,L_1)$ acquire additional integer indices $i^r(x) \in \Z$. In the case of \eqref{eq:hom01} where the two underlying curves agree, the difference of their $1/r$-gradings is constant, $a_1(x) - a_0(x) = \pi d$ for some $d \in \Z$. One then sets $i^r(e) = d$, $i^r(q) = d+r$. In the other situation \eqref{eq:cf}, let $\alpha \in (0,\pi)$ be the angle with which our Lagrangian submanifolds meet at $x$, counted clockwise from $\mathit{TL}_{0,x}$ to $\mathit{TL}_{1,x}$, and let $a_0(x),a_1(x)$ be the $1/r$-gradings at that point. Define
\begin{equation}
i^r(x) = \frac{r\alpha(x)+a_1(x)-a_0(x)}{\pi}.
\end{equation}
If $r$ is odd, which is when orientations actually matter in \eqref{eq:phase}, the parity of $i^r$ agrees with the previously used $\Z/2$ grading.

Consider a moduli space $\scrM(x_0,\dots,x_d)$ which enters into the definition of the $A_\infty$-structure of $\F(M,D)$. For simplicity, we assume that the Lagrangian submanifolds $(L_0,\dots,L_d)$ involved are pairwise distinct (similar arguments would apply in the more general case). Take a point $(\zeta^0,\dots,\zeta^d,u)$ in our moduli space. Since the boundary of $u$ lies in $M \setminus D$, for any $z \in D$ there is a well-defined degree $\deg(u,z)$, namely the multiplicity with which $u$ hits $z$. It is nonnegative, and vanishes iff $u^{-1}(z) = \emptyset$. Using the fact that our point has virtual dimension zero, and the index formula \cite[Section 11]{seidel04}, one sees that
\begin{equation} \label{eq:r-index}
i^r(x_0) - i^r(x_1) - \cdots - i^r(x_d) = r(2-d) + 2\sum_{z \in D} \mathrm{ord}(\eta^r,z) \deg(u,z).
\end{equation}

Suppose for simplicity that $\eta$ vanishes to the same order $m = \mathrm{ord}(\eta^r,z) > 0$ at every point $z \in D$. One way to encode \eqref{eq:r-index} is to equip $\mathit{CF}^*(L_0,L_1)$ with the grading given by the $i^r$, and then to write the composition maps in $\F(M,D)$ as
\begin{equation} \label{eq:decompose-mu}
\mu^i = \mu^i_0 + \mu^i_1 + \cdots
\end{equation}
where the subscript denotes the total degree of $u$ over $D$, and the term $\mu^i_k$ has degree $r(2-i) + 2km$ with respect to the $i^r(x_k)$. Equivalently, one could add a formal variable $\hbar$ of degree $2m$ to form $\mathit{CF}^*(L_0,L_1)[[\hbar]]$, and then get a homogeneous map $\mu^i_0 + \hbar \mu^i_1 + \cdots$ of degree $r(2-i)$. The sum in \eqref{eq:decompose-mu} is of course finite when applied to any given $(L_0,\dots,L_d)$. Correspondingly, each graded piece of $\mathit{CF}^*(L_0,L_1)[[\hbar]]$ only contains finitely many powers of $\hbar$.

\section{Orbifolds}

Part of our argument will involve working equivariantly with respect to the action of a finite group on the target surface. Even though it makes no fundamental difference, it can be more intuitive to think of this as working on the orbifold quotient, so we'll give a short discussion which takes this point of view into account. Take $M$ as before. Let $\Gamma$ be a finite group acting effectively on $M$ by orientation-preserving diffeomorphisms, and let $D$ be the finite set of points where the action is not free. Write $\bar{M} = M/\Gamma$ for the orbifold quotient, and $\bar{D}$ for its finite set of orbifold points. We can choose $\omega$ and $\theta$ to be $\Gamma$-invariant.
%\end{quote}
%\begin{figure}[hb]
%\begin{centering}
%\epsfig{file=tear.eps}
%\caption{\label{fig:teardrop}}
%\end{centering}
%\end{figure}
%Moreover, this has the property that there are no teardrops in $\bar{M} \setminus \bar{D}$ with boundary on $\bar{L}$ (see Figure %\ref{fig:teardrop} for a visual explanation of the terminology).

Take an embedded balanced curve $L \subset M \setminus D$, such that the collection $\{\gamma(L) \;:\; \gamma \in \Gamma\}$ is in general position, hence can be assumed to be part of our family $\scrL$. Then, the image curve $\bar{L} \subset \bar{M} \setminus \bar{D}$ is immersed and has only transversal double points. Equip $L$ with an orientation, a {\em Spin} structure, and a Morse function $f$ as in Section \ref{sec:technical}. We then define
\begin{equation} \label{eq:orbifold-floer-cochain}
\mathit{CF}^*(\bar{L},\bar{L}) = \mathit{CM}^*(f) \oplus \bigoplus_{\gamma \neq 1} \mathit{CF}^*(L,\gamma(L)).
\end{equation}
Denote the summands by $\mathit{CF}^*(\bar{L},\bar{L})^\gamma$, where the first term corresponds to $\gamma = 1$. We have formulated the definition in terms of $L$ for convenience. In terms of $\bar{L}$ itself, we have two generators corresponding to the classical Morse chain complex, as well as a pair of generators for each selfintersection point (this is a general feature of Floer cohomology for immersed Lagrangian submanifolds, see \cite{akaho05,akaho-joyce08}). Suppose that $(\bar{x}_0,\dots,\bar{x}_d)$ are generators, belonging to the $(\gamma_0,\dots,\gamma_d)$ summands of \eqref{eq:orbifold-floer-cochain}, and $(x_0,\dots,x_d)$ their obvious lifts to $M$. The associated moduli space $\scrM(\bar{x}_0,\dots,\bar{x}_d)$ is empty if $\gamma_0 \neq \gamma_1 \gamma_2 \cdots \gamma_d$, and can be identified with $\scrM(x_0,x_1,\gamma_1(x_2),\gamma_1\gamma_2(x_3),\dots,\gamma_1\cdots\gamma_{d-1}(x_d))$
in the remaining case. Here, the assumption is that the Riemann surface structure, and other auxiliary choices, are made equivariantly. In simple situations where points of that moduli space are holomorphic discs, the image of any such disc under the quotient map $M \rightarrow \bar{M}$ is an ``orbifold holomorphic disc'' (meaning that it has appropriate ramification at all points of $\bar{D}$), and conversely all orbifold holomorphic discs lift to appropriate objects in $M$.

\begin{remark}
Transversality may seem to be an issue in the equivariant context, but a quick reflection shows that this is not the case in the present context, where the group action is free on the set of objects $\{\gamma(L)\}$. If one uses virtual perturbation methods as in \cite{fooo}, multivalued perturbations are built into the framework, and that naturally allows equivariance with respect to any given finite group action. Alternatively, let's consider using explicit perturbations as in \cite{seidel04}. Such perturbations are given by inhomogeneous terms which vary on the domain of our holomorphic maps, which means that we are looking at equations of the class
\begin{equation}
\bar\partial u = \nu(z,u(z)).
\end{equation}
The group $\Gamma$ acts on the target space, but leaves $z$ invariant, and that allows enough freedom even if one takes $\nu$ to be $\Gamma$-invariant and to vanish near $D$. This is also easy to see if one thinks in terms of $\bar{M}$.
\end{remark}

It should be clear from the definition that there is a simple relationship between the $A_\infty$-structures obtained by looking at $\bar{L} \subset \bar{M}$ and at all preimages $\gamma(L) \subset M$. To state this in a simple algebraic way, suppose that $\Gamma$ is abelian, and consider its character group $G = \mathrm{Hom}(\Gamma,\C^*)$. We then have an action of $G$ on $\mathit{CF}^*(\bar L,\bar L)$, defined by $g \cdot x = g(\gamma) x$ for $x \in \mathit{CF}^*(\bar L,\bar L)^\gamma$. By construction, the $A_\infty$-structure is equivariant with respect to this action, and we have
\begin{equation} \label{eq:rtimes}
\bigoplus_{\gamma_0,\gamma_1} \mathit{CF}^*(\gamma_0(L),\gamma_1(L)) = \mathit{CF}^*(\bar{L},\bar{L}) \rtimes G,
\end{equation}
compatibly with the $A_\infty$-structures, where $\rtimes$ is the semidirect product. This is familiar from other instances of mirror symmetry, see for instance \cite{seidel03b}. Conversely, the left hand side of \eqref{eq:rtimes} carries an action of $\Gamma$, whose invariant part can be identified with $\mathit{CF}^*(\bar{L},\bar{L})$.

\begin{remark}
Suppose that $\bar{M}$ has genus zero. Then $\Gamma = H_1^{\mathrm{orb}}(\bar{M})$ is a finite abelian group, and one can write $\bar{M} = M/\Gamma$. The associated group action can be interpreted in more familiar, if somewhat more abstract, terms as follows. Suppose that we introduced a larger version of the Fukaya category of the orbifold, in which our curves carry flat $\C^*$-bundles. Then the Picard groupoid, which consists of flat $\C^*$-bundles defined on the entire orbifold, acts on that category by tensor product. In particular, if we have a curve $\bar{L}$ which is fixed under this action in a suitable sense, which means that its class in $H_1^{\mathrm{orb}}(\bar{M})$ is zero, then $G = H^1_{\mathrm{orb}}(\bar{M};\C^*)$ acts on its endomorphism space $\mathit{CF}^*(\bar{L},\bar{L})$.
\end{remark}

\section{The genus two case\label{sec:g2}}

From now on, we return to the specific case where $M$ has genus $2$. Initially, it will be convenient to represent this as a double cover of $S^2 = \C \cup \{\infty\}$ branched over six points, which are the fifth roots of unity and $0$. A nonseparating simple closed curve in $M$ which is invariant under the hyperelliptic involution projects to an embedded path in $S^2$ which connects two of the branch points. Similarly, the Dehn twist along the curve projects to the half-twist along the associated path. We start with the configuration of five curves in $M$ whose images in $S^2$ form the pentagram, Figure \ref{fig:kabbalah}(i).
\begin{figure}[hb]
\begin{centering}
\begin{picture}(0,0)%
\includegraphics{kabbalah.pstex}%
\end{picture}%
\setlength{\unitlength}{3552sp}%
\begingroup\makeatletter\ifx\SetFigFont\undefined%
\gdef\SetFigFont#1#2#3#4#5{%
  \reset@font\fontsize{#1}{#2pt}%
  \fontfamily{#3}\fontseries{#4}\fontshape{#5}%
  \selectfont}%
\fi\endgroup%
\begin{picture}(5163,2182)(1036,-1634)
\put(5596,-811){\makebox(0,0)[lb]{\smash{{\SetFigFont{11}{13.2}{\rmdefault}{\mddefault}{\updefault}{\color[rgb]{0,0,0}$K_4$}%
}}}}
\put(3901,389){\makebox(0,0)[lb]{\smash{{\SetFigFont{11}{13.2}{\rmdefault}{\mddefault}{\updefault}{\color[rgb]{0,0,0}(ii)}%
}}}}
\put(1051,389){\makebox(0,0)[lb]{\smash{{\SetFigFont{11}{13.2}{\rmdefault}{\mddefault}{\updefault}{\color[rgb]{0,0,0}(i)}%
}}}}
\put(1296,-566){\makebox(0,0)[lb]{\smash{{\SetFigFont{11}{13.2}{\rmdefault}{\mddefault}{\updefault}{\color[rgb]{0,0,0}$L_1$}%
}}}}
\put(2776,-961){\makebox(0,0)[lb]{\smash{{\SetFigFont{11}{13.2}{\rmdefault}{\mddefault}{\updefault}{\color[rgb]{0,0,0}$L_3$}%
}}}}
\put(2656,-136){\makebox(0,0)[lb]{\smash{{\SetFigFont{11}{13.2}{\rmdefault}{\mddefault}{\updefault}{\color[rgb]{0,0,0}$L_4$}%
}}}}
\put(2926,-586){\makebox(0,0)[lb]{\smash{{\SetFigFont{11}{13.2}{\rmdefault}{\mddefault}{\updefault}{\color[rgb]{0,0,0}$L_2$}%
}}}}
\put(1466,-961){\makebox(0,0)[lb]{\smash{{\SetFigFont{11}{13.2}{\rmdefault}{\mddefault}{\updefault}{\color[rgb]{0,0,0}$L_5$}%
}}}}
\put(4801, 14){\makebox(0,0)[lb]{\smash{{\SetFigFont{11}{13.2}{\rmdefault}{\mddefault}{\updefault}{\color[rgb]{0,0,0}$K_1$}%
}}}}
\put(4501,-811){\makebox(0,0)[lb]{\smash{{\SetFigFont{11}{13.2}{\rmdefault}{\mddefault}{\updefault}{\color[rgb]{0,0,0}$K_2$}%
}}}}
\put(5026,-1561){\makebox(0,0)[lb]{\smash{{\SetFigFont{11}{13.2}{\rmdefault}{\mddefault}{\updefault}{\color[rgb]{0,0,0}$K_3$}%
}}}}
\end{picture}%
\caption{\label{fig:kabbalah}}
\end{centering}
\end{figure}

\begin{lemma} \label{th:l-split-generate}
The curves $\{L_1,\dots,L_5\}$, equipped with nontrivial ${\mathit Spin}$ structures, split-generate $D^\pi\F(M)$.
\end{lemma}

\proof Consider first the collection $\{K_1,K_2,K_3,K_4\}$ from Figure \ref{fig:kabbalah}(ii). The Dehn twists $\tau_{K_i}$ define a homomorphism from the braid group $\mathit{Br}_5$ to the mapping class group of $M$. Its kernel is infinite cyclic and generated by the central element $\Delta^4 \in \mathit{Br}_5$, which in particular means that its image $(\tau_{K_4}\cdots\tau_{K_1})^{10}$ is isotopic to the identity \cite{birman72}. Lemma \ref{th:karoubi-2} shows that the $K_i$ split-generate $D^\pi\F(M)$. On the other hand, one easily checks by hand that $\tau_{L_5}\cdots\tau_{L_1}(K_2)$ is isotopic to $K_2[1]$, and analogously for the other $K_i$. By Lemma \ref{th:karoubi-1}, each $K_i$ is split-generated by $\{L_1,\dots,L_5\}$, which completes the argument. \qed

For the main computation we switch to a different picture. Take the action of $\Sigma = \Z/5$ on $M$ which projects to the rotational action on $S^2$. The orbifold quotient $\bar{M} = M/\Sigma$ is a sphere with three orbifold points $\bar{D}$. Each $L = L_i$ projects to the same immersed curve $\bar{L} \subset \bar{M}$. We will use the previously introduced techniques to partially compute the $A_\infty$-structure on $\mathit{CF}^*(\bar{L},\bar{L})$. Note that all such computations can be equivalently thought of as being carried out on $M$, which means that the orbifold structure does not really introduce any new technical issues.
\begin{figure}
\begin{centering}
\begin{picture}(0,0)%
\includegraphics{generators.pstex}%
\end{picture}%
\setlength{\unitlength}{3355sp}%
\begingroup\makeatletter\ifx\SetFigFont\undefined%
\gdef\SetFigFont#1#2#3#4#5{%
  \reset@font\fontsize{#1}{#2pt}%
  \fontfamily{#3}\fontseries{#4}\fontshape{#5}%
  \selectfont}%
\fi\endgroup%
\begin{picture}(5381,3774)(889,-3223)
\put(1951,-1186){\makebox(0,0)[lb]{\smash{{\SetFigFont{10}{12.0}{\familydefault}{\mddefault}{\updefault}{\color[rgb]{0,0,0}$x_1,\,\bar{x}_1$}%
}}}}
\put(4651,-1186){\makebox(0,0)[lb]{\smash{{\SetFigFont{10}{12.0}{\familydefault}{\mddefault}{\updefault}{\color[rgb]{0,0,0}$x_3,\,\bar{x}_3$}%
}}}}
\put(3226,-2836){\makebox(0,0)[lb]{\smash{{\SetFigFont{10}{12.0}{\familydefault}{\mddefault}{\updefault}{\color[rgb]{0,0,0}$x_2,\,\bar{x}_2$}%
}}}}
\put(3976,-1336){\makebox(0,0)[lb]{\smash{{\SetFigFont{10}{12.0}{\rmdefault}{\mddefault}{\updefault}{\color[rgb]{0,0,0}$q$}%
}}}}
\put(3676,-1711){\makebox(0,0)[lb]{\smash{{\SetFigFont{10}{12.0}{\rmdefault}{\mddefault}{\updefault}{\color[rgb]{0,0,0}$e$}%
}}}}
\end{picture}%
\caption{\label{fig:generators}}
\end{centering}
\end{figure}

{\em Generators:} We have the Morse-theoretical generators $e$ (minimum, even) and $q$ (maximum, odd), together with a pair of generators coming from each self-intersection point, which we denote by $\bar{x}_k$ (even) and $x_k$ (odd), see Figure \ref{fig:generators} (we have perturbed that picture slightly to make the self-intersections more visible; the more natural picture would be the one with full symmetry between the front and back faces). Take $\Gamma = H_1^{\mathrm{orb}}(\bar{M})$, which one thinks of as the quotient of $(\Z/5)^3$ by its diagonal subgroup $\Z/5$. Since the class of our immersed curve in $\Gamma$ is trivial, the generators of the Floer cochain complex come labeled by weights which are elements of $\Gamma$. Moreover, there is a nontrivial holomorphic section $\eta^3$ of $(T^*\!\bar{M})^{\otimes 3}$, unique up to nonzero scalars, which has a double zero at each point of $\bar{D}$ (this is in fact the same as a meromorphic section of $(T^*\!S^2)^{\otimes 3}$ with double poles at our three points, with the different order due to considering it in orbifold charts: if $z = w^5$, then $z^{-2}dz^3 = 5^3 w^2 dw^3$). As a consequence, the generators acquire additional integer indices. All this data can be listed as follows:
\begin{equation}
\begin{aligned}
&
\begin{array}{l|l|l|l|l}
\text{generator} & e & x_1 & x_2 & x_3 \\ \hline
\text{weight} & (0,0,0) & (1,0,0) & (0,1,0) & (0,0,1) \\
\text{index} & 0 & 1 & 1 & 1
\end{array} \\
&
\begin{array}{l|l|l|l|l}
\text{generator} & \bar{x}_1 & \bar{x}_2 & \bar{x}_3 & q \\ \hline
\text{weight} & (0,1,1) & (1,0,1) & (1,1,0) & (1,1,1) \\
& = (-1,0,0) & = (0,-1,0) & = (0,0,-1) & = (0,0,0) \\
\text{index} & 2 & 2 & 2 & 3
\end{array}
\end{aligned}
\end{equation}
We know that the $A_\infty$-structure is homogeneous with respect to weights. This immediately implies that $\mu^1 = 0$. Moreover, we have a decomposition as in \eqref{eq:decompose-mu}, where $\mu^i_k$ has degree $6 - 3i + 4k$ with respect to the indices above. Concretely, $\mu^i_0$ counts the contributions from polygons in $\bar{M} \setminus \bar{D}$; the next term $\mu^i_1$, that from polygons which meet $\bar{D}$ exactly once and have fivefold ramification at that point (which is the minimal order prescribed by the orbifold structure).

{\em Triangles:} For degree reasons, $\mu^2_k = 0$ for all $k>0$. According to \eqref{eq:e2} and \eqref{eq:e3}, the contributions of the constant triangles are
\begin{equation}
\begin{aligned}
& \mu^2(x_i,e) = x_i = -\mu^2(e,x_i), \\ & \mu^2(\bar{x}_i,e) = \bar{x}_i = \mu^2(e,\bar{x}_i), \\ & \mu^2(q,e) = q = -\mu^2(e,q), \\
& \mu^2(q,q) = 0, \\
& \mu^2(x_i,\bar{x}_i) = q = - \mu^2(\bar{x}_i,x_i).
\end{aligned}
\end{equation}
There are six (if one counts the ordering of their corners; otherwise, only two) non-constant triangles avoiding $\bar{D}$. To determine the sign of their contribution, we need to choose generic points $\circ$ on $L$, which represent the nontrivial $\mathit{Spin}$ structure. Those being as in Figure \ref{fig:generators}, we get
\begin{equation}
\begin{aligned}
& \mu^2_0(x_1,x_2) = \bar{x}_3 = -\mu^2_0(x_2,x_1), \\
& \mu^2_0(x_2,x_3) = \bar{x}_1 = -\mu^2_0(x_3,x_2), \\
& \mu^2_0(x_3,x_1) = \bar{x}_2 = -\mu^2_0(x_1,x_3).
\end{aligned}
\end{equation}
The triangle on the front part of Figure \ref{fig:generators} goes through the $e$, hence can be thought of as a holomorphic map from the four-punctured disc which is smooth at one of the marked points. The resulting contribution is
\begin{equation}
\mu^3_0(x_3,x_2,x_1) = -e.
\end{equation}
All other expressions $\mu^3_0(x_{i_3},x_{i_2},x_{i_1})$ are zero: any such product can only be a multiple of $e$, for degree reasons, but the relevant spaces $\scrM(e,x_{i_3},x_{i_2},x_{i_1})$ are empty.
\begin{figure}
\begin{centering}
\begin{picture}(0,0)%
\includegraphics{1-penta.pstex}%
\end{picture}%
\setlength{\unitlength}{3552sp}%
\begingroup\makeatletter\ifx\SetFigFont\undefined%
\gdef\SetFigFont#1#2#3#4#5{%
  \reset@font\fontsize{#1}{#2pt}%
  \fontfamily{#3}\fontseries{#4}\fontshape{#5}%
  \selectfont}%
\fi\endgroup%
\begin{picture}(2488,1749)(889,-1198)
\put(2101,-593){\makebox(0,0)[lb]{\smash{{\SetFigFont{11}{13.2}{\familydefault}{\mddefault}{\updefault}{\color[rgb]{0,0,0}}%
}}}}
\end{picture}%
\caption{\label{fig:1-penta}}
\end{centering}
\end{figure}
\begin{figure}
\begin{centering}
\begin{picture}(0,0)%
\includegraphics{2-penta.pstex}%
\end{picture}%
\setlength{\unitlength}{3552sp}%
\begingroup\makeatletter\ifx\SetFigFont\undefined%
\gdef\SetFigFont#1#2#3#4#5{%
  \reset@font\fontsize{#1}{#2pt}%
  \fontfamily{#3}\fontseries{#4}\fontshape{#5}%
  \selectfont}%
\fi\endgroup%
\begin{picture}(4456,4338)(4109,-3752)
\put(5191,-387){\makebox(0,0)[lb]{\smash{{\SetFigFont{11}{13.2}{\familydefault}{\mddefault}{\updefault}{\color[rgb]{0,0,0}2}%
}}}}
\put(4790,-270){\makebox(0,0)[lb]{\smash{{\SetFigFont{11}{13.2}{\familydefault}{\mddefault}{\updefault}{\color[rgb]{0,0,0}1}%
}}}}
\put(4763,-2150){\makebox(0,0)[lb]{\smash{{\SetFigFont{11}{13.2}{\familydefault}{\mddefault}{\updefault}{\color[rgb]{0,0,0}1}%
}}}}
\put(5206,-2206){\makebox(0,0)[lb]{\smash{{\SetFigFont{11}{13.2}{\familydefault}{\mddefault}{\updefault}{\color[rgb]{0,0,0}1}%
}}}}
\end{picture}%
\caption{\label{fig:2-penta}}
\end{centering}
\end{figure}

{\em Pentagons:} There are six ``pentagons'' hidden in the picture. Each of them hits exactly one of the points of $\bar{D}$ and has fivefold ramification there, and no ramification elsewhere, which means that it lifts to a genuine immersed pentagon in $M$. Figure \ref{fig:1-penta} shows the image of one of the pentagons in the orbifold picture. Figure \ref{fig:2-penta} shows two pentagons lifted to $M$ and then projected to $S^2$ under the hyperelliptic quotient (the dots in the left hand picture are the corners, and the numbers local degrees; the right hand side shows a smoothed version of the boundary curve). We are not interested in their effect on $\mu^4_1$ (this happens to cancel, but its vanishing is in fact a consequence of the previous computations and the $A_\infty$ structure equations, hence yields no new information). Instead, we take the three pentagons whose boundary goes through $e$, and determine their contributions to $\mu^5_1$, namely:
\begin{equation}
\mu^5_1(x_1,x_1,x_1,x_1,x_1) = \mu^5_1(x_3,x_3,x_3,x_3,x_3) = -e, \; \mu^5_1(x_2,x_2,x_2,x_2,x_2) = e.
\end{equation}
The other $\mu^5_1(x_{i_5},x_{i_4},x_{i_3},x_{i_2},x_{i_1})$, where the indices are not all the same, vanish.

Identify $\mathit{CF}^*(\bar{L},\bar{L}) \iso \Lambda(V)$ by mapping the generators as follows:
\[
\begin{array}{l|l|l|l|l|l|l|l|l}
\mathit{CF}^*(\bar{L},\bar{L}) & e & x_1 & x_2 & x_3 & \bar{x}_1 & \bar{x}_2 & \bar{x}_3 & q \\
\hline
\Lambda(V) & 1 & -\xi_1 & \xi_2 & -\xi_3 & \xi_2 \wedge \xi_3 & \xi_1 \wedge \xi_3 & \xi_1 \wedge \xi_2 & -\xi_1 \wedge \xi_2 \wedge \xi_3
\end{array}
\]
Then, the data above fit precisely into Proposition \ref{th:classification}, proving that $\mathit{CF}^*(\bar{L},\bar{L})$ is $A_\infty$-isomorphic to $\A$. Moreover, this isomorphism is equivariant with respect to the action of $G = \mathrm{Hom}(\Gamma,\C^*) \subset \mathit{SL}(V)$. The orbifold covering $M \rightarrow \bar{M}$ is classified by a surjective homomorphism $\Gamma \rightarrow \Sigma$, whose dual is the subgroup $Z \subset G$. In view of \eqref{eq:rtimes} we then inherit an $A_\infty$-isomorphism
\begin{equation}
\bigoplus_{i,j} \mathit{CF}^*(L_i,L_j) \iso \A \rtimes Z.
\end{equation}
We can consider the left hand as an $A_\infty$-category with a single object, which is the formal direct sum of the $L_i$. Using Lemma \ref{th:l-split-generate} we arrive at the following description of the Fukaya category $M$ on the derived level:

\begin{corollary} \label{th:describes-f}
$D^\pi\F(M) \iso D^\pi(\A \rtimes Z)$.
\end{corollary}

\section{Koszul duality\label{sec:koszul}}

In this section, we again allow $V = \C^n$ for any $n$. Consider the space of differential forms $\Omega(V) = \C[V] \otimes \Lambda(V^\vee)$, with the grading reversed (negative), and equip it with the differential $\iota_\eta$ given by contraction with the Euler vector field $\eta = \sum_i v_i \xi_i$. Consider the dga
\begin{equation} \label{eq:b-dga}
B = \mathrm{Hom}_{\C[V]}(\Omega(V),\Omega(V))
\end{equation}
with the induced differential $\partial$. Since $\Omega(V)$ is just the standard free Koszul resolution of the simple $\C[V]$-module $\C$, we have $H(B) = \mathrm{Ext}^*_{\C[V]}(\C,\C) \iso \Lambda(V)$. To write things down on the cochain level, let's identify
\begin{equation}
B = \mathrm{Hom}_\C(\Lambda(V^\vee),\Omega(V)) = \Omega(V)  \otimes \Lambda(V),
\end{equation}
and correspondingly write the differential as
\begin{equation}
\begin{aligned}
& \partial: \Omega(V) \otimes \Lambda(V) \longrightarrow \Omega(V) \otimes \Lambda(V), \\
& \partial(f\beta \otimes \theta) = \sum_k v_kf \, \iota_{\xi_k}\beta \otimes \theta + (-1)^{|\beta|-1}
v_k f\beta \otimes \xi_k \wedge \theta.
\end{aligned}
\end{equation}
The projection $p: \Omega(V) \otimes \Lambda(V) \rightarrow \Lambda(V)$, where the target $\Lambda(V)$ carries the zero differential, is a chain homomorphism. In converse direction we have the inclusion
\begin{equation}
\begin{aligned}
& i: \Lambda(V) \longrightarrow \Omega(V) \otimes \Lambda(V), \\
& i(\theta) = \sum_{p \geq 0} \sum_{j_1 < \dots < j_p} dv_{j_1} \wedge \cdots \wedge dv_{j_p} \otimes
\xi_{j_p} \wedge \cdots \wedge \xi_{j_1} \wedge \theta.
\end{aligned}
\end{equation}
Finally, there is a chain homotopy $h$ between the identity and $p \, i$:
\begin{equation}
\begin{aligned}
& h: \Omega(V) \otimes \Lambda(V) \longrightarrow \Omega(V) \otimes \Lambda(V), \\
& h(f \beta \otimes \theta) = 0 \;\; \text{if $f\beta$ is a multiple of $1$, and otherwise} \\
& = \sum_{p \geq 0}
{\textstyle \frac{p!}{w(w+1)\cdots (w+p)}}
\sum_{j_1 < \dots < j_p }
df \wedge \beta \wedge dv_{j_1} \wedge \cdots \wedge dv_{j_p} \otimes \xi_{j_p} \wedge \cdots \wedge \xi_{j_1} \wedge \theta
\end{aligned}
\end{equation}
where $w = r+s$ for $f \in \Sym^r(V^\vee)$ and $\beta \in \Lambda^s(V^\vee)$ ($w$ is the weight of $f\beta$ with respect to the diagonal $\C^*$-action on $V^\vee$). Moreover, it is tacitly assumed that $h(f\beta \otimes \gamma)$ vanishes if $w = 0$. The entire structure constructed above also satisfies the so-called side conditions
\begin{equation} \label{eq:side-condition}
h^2 = 0, \;\; p \, h = 0, \;\; h \, i = 0.
\end{equation}

Returning to the original definition \eqref{eq:b-dga}, one could also write $i(\theta)$ as the action of $\theta$ by contraction, $\iota_\theta: \Omega(V) \rightarrow \Omega(V)$. From this, it follows that $i$ is a map of differential graded algebras. Hence $B$ is formal. This is the prototypical instance of a more general homogeneity phenomenon, known as Koszul duality (see for instance \cite{beilinson-ginzburg-soergel96}). We will be interested in a deformed version of this statement. Namely, take a one-form $\gamma = \sum_k g_k dv_k \in \Omega^1(V)$, and change the differential on $\Omega(V)$ to $\iota_\eta - \gamma \wedge \cdot$, which of course reduces the grading to $\Z/2$. The square of the new differential is multiplication with the function
\begin{equation} \label{eq:w-as-eta}
W = -\gamma(\eta) \in \C[V],
\end{equation}
which is a central element. Hence, the induced differential on $B$ does indeed square to zero. We denote this differential by $\tilde\partial$, and the resulting dga structure by $\B$. Explicitly,
\begin{equation}
(\tilde\partial - \partial) (f\beta \otimes \theta) = -f\, \gamma \wedge \beta \otimes \theta + (-1)^{|\beta|-1}
\sum_k g_k f \beta \otimes \iota_{dv_k} \theta.
\end{equation}

Starting from this, the Homological Perturbation Lemma \cite{markl04} constructs an induced $\Z/2$-graded $A_\infty$-structure $\A$ on $A = \Lambda(V)$, together with an $A_\infty$-quasi-isomorphism from that structure to $\B$. An explicit formula for the differential is
\begin{equation} \label{eq:chain-mu1}
\mu^1(a) = (-1)^{|a|} p (\tilde\partial-\partial) i(a) + (-1)^{|a|} p(\tilde\partial-\partial)h(\tilde\partial-\partial) i(a) + \cdots
\end{equation}
The general formula for $\mu^d$ is as a sum over ribbon trees with a root and $d$ leaves, whose vertices may have valencies $2$ or $3$. Take such a tree and orient it in a way pointing from the leaves to the root, then attach an operation to each vertex and edge, as follows:
\begin{equation}\left\{
\begin{aligned}
& \text{for a bivalent vertex,} && b \mapsto (-1)^{|b|} (\tilde\partial-\partial)(b): B \rightarrow B, \\
& \text{for a trivalent vertex,} && (b_2,b_1) \mapsto (-1)^{|b_1|} b_2b_1: B \otimes B \rightarrow B, \\
& \text{for a finite edge,} && b \mapsto (-1)^{|b|-1} h(b): B \rightarrow B, \\
& \text{for a semi-infinite incoming edge,}, && a \mapsto i(a): A \rightarrow B, \\
& \text{for a semi-infinite outgoing edge,} && b \mapsto p(b): B \rightarrow A.
\end{aligned}\right.
\end{equation}
Then compose these operations as prescribed by the tree itself, to get a multilinear map $A^{\otimes d} \rightarrow A$. For instance, the terms in \eqref{eq:chain-mu1} arise from the linear trees (chains of bivalent vertices, with one semi-infinite incoming and another semi-infinite outgoing end). Because both $\tilde\partial - \partial$ and $h$ decrease the grading, trees containing sufficiently long chains contribute zero, and therefore all resulting sums are finite.

\begin{lemma} \label{th:low-order}
Fix some $r \geq 0$. Suppose that all $g_k$ lie in $F_{r-1}\C[V]$, which means that they do not contain monomials of order $< r-1$. Then the resulting $A_\infty$-structure on $A$ agrees with the trivial (formal) one up to order $r-1$.
\end{lemma}

To see this, consider the grading of $B$ by the order of its symmetric algebra part. This grading is preserved by the product structure, decreased by one under $h$, and increased by at least $r-1$ under $\tilde\partial-\partial$. Hence, the multilinear map arising from a tree with $d$ leaves and $k$ bivalent vertices can be nonzero only if $k \leq (d-2)/(r-2)$. On the other hand, trees with $d \geq 3$ and $k = 0$ contribute zero, because $h(i(a_2)i(a_1)) = h(i(a_2a_1)) = 0$.

\begin{remark}
The construction above gives explicit formulae for the entire $A_\infty$-deformation of the exterior algebra induced by the superpotential $W$. These could be useful in other situations, such as ones considered in \cite{aspinwall-katz06}. It should also be mentioned that there is another possible way of obtaining such formulae, namely by applying \eqref{eq:map-alpha} and Kontsevich formality. I have not investigated the relation between the two approaches.
\end{remark}

We now return to the usual special case, where $V = \C^3$ and $W$ is as in \eqref{eq:w}. Take
\begin{equation}
g_1 = -v_2 v_3/3 + v_1^4, \;\; g_2 = -v_1v_3/3 + v_2^4, \;\; g_3 = -v_1v_2/3 + v_3^4.
\end{equation}

\begin{proposition} \label{th:endomorphism-characterized}
The resulting $A_\infty$-algebra $\A$, with the obvious action of $G \subset \mathit{GL}(V)$, belongs to the quasi-isomorphism class singled out in Proposition \ref{th:classification}.
\end{proposition}

For this, it is convenient to make some temporary changes. Let's first modify the grading of $B$ by giving the summand $\Sym^i(V^\vee) \otimes \Lambda^j(V^\vee) \otimes \Lambda^k(V)$ degree $2i-j+k$. The product is still compatible with this grading, but $\partial$ has degree $3$ and $h$ has degree $-3$. To make the remaining term $\tilde\partial-\partial$ have degree $3$ as well, introduce a formal parameter $\hbar$ of degree $-4$, and write $g_1 = -v_2 v_3/3 + \hbar v_1^4$ and similarly for the other $g_k$. Then, the resulting $A_\infty$-operations can be written as sums of terms $\mu^d_k$, coming with $\hbar^k$, of degrees $6-3d+4k$. Moreover, since $\tilde\partial-\partial$ is $G$-equivariant, and all the other data are equivariant for the entire group $\mathit{GL}(V)$, the $A_\infty$-operations inherit $G$-symmetry.

Lemma \ref{th:low-order} shows that $\mu^1$ vanishes and $\mu^2$ is the standard wedge product. It remains to compute the two higher order compositions in \eqref{eq:higher-35}. This is elementary, using the explicit tree summation formulae discussed above, and the result is precisely as required. We omit the details, referring instead to \cite{seidel08x}.

\section{Matrix factorizations\label{sec:matrix}}

For the duration of this section, we consider the more general case where $V = \C^n$ for any $n$, and $W \in \C[V]$ is a polynomial such that the hypersurface $W^{-1}(0)$ has a single singular point, which lies at the origin. Orlov \cite{orlov04} associates to this hypersurface the category
\begin{equation} \label{eq:db-sing}
D^b_{\mathrm{sing}}(W^{-1}(0)) = D^b(W^{-1}(0))/\mathrm{Perf}(W^{-1}(0)).
\end{equation}
The notation here is that $D^b(W^{-1}(0))$ is the bounded derived category of coherent sheaves; $\mathrm{Perf}(W^{-1}(0))$ the full triangulated subcategory of perfect complexes; and the quotient is localization with respect to the class of morphisms whose cones lie in $\mathrm{Perf}(W^{-1}(0))$. By definition $D^b_{\mathrm{sing}}(W^{-1}(0))$ is triangulated. A deeper fact, based on cohomological properties of hypersurfaces, is that it is $\Z/2$-graded, which as before means that the twofold shift is isomorphic to the identity. Finally, the categories $D^b_{\mathrm{sing}}(W^{-1}(0))$ are not split-closed in general, but we can take the split-closure $D^\pi_{\mathrm{sing}}(W^{-1}(0))$, which is again naturally triangulated \cite{balmer-schlichting01}.

\begin{lemma} \label{th:skyscraper-generate}
$D^\pi_{\mathrm{sing}}(W^{-1}(0))$ is split-generated by the skyscraper sheaf at the origin, $\EuScript{S}_{W^{-1}(0),0}$.
\end{lemma}

\proof A result of Orlov \cite{orlov08b} says that any object in that category is a direct summand of the image of an object of $D^b_{\mathrm{sing}}(W^{-1}(0))$ whose cohomology sheaves are supported on the singular locus of $W^{-1}(0)$, in this case the origin. On the other hand, any such complex can be built from shifted copies of $\EuScript{S}_{W^{-1}(0),0}$ through repeated mapping cones. Since the projection functor $D^b(W^{-1}(0)) \rightarrow D^b_{\mathrm{sing}}(W^{-1}(0))$ is exact by definition, this behaviour transfers to the quotient category. \qed

Matrix factorizations \cite{eisenbud80}, which historically predate $D^b_{\mathrm{sing}}(W^{-1}(0))$, can be used to construct a chain level model for that category. With $W$ as before, a matrix factorization is a $\Z/2$-graded projective $\C[V]$-module $E$ together with an odd $\C[V]$-linear differential $\delta_E$ such that $\delta_E^2 = W \cdot \mathrm{id}_E$. Matrix factorizations form a $\Z/2$-graded differential graded category $\mathit{MF}(W)$, and this admits mapping cones, hence the cohomological category $H^0(\mathit{MF}(W))$ is naturally triangulated.

\begin{theorem}[\protect{Orlov \cite[Theorem 3.9]{orlov04}}] \label{th:orlov}
There is an equivalence of triangulated categories,
\begin{equation} \label{eq:orlov-equiv}
H^0(\mathit{MF}(W)) \iso D^b_{\mathrm{sing}}(W^{-1}(0)).
\end{equation}
On the level of objects, this takes a matrix factorization $E$ to the coherent sheaf corresponding to the $\C[V]/W$-module $\mathrm{coker}(\delta_E^1: E^1 \rightarrow E^0)$.
\end{theorem}

We will now make the connection with the material from the previous section. Suppose that $W$ is written in the form \eqref{eq:w-as-eta}, and take $E = \Omega(V)$ with its natural $\Z/2$-grading, and with the differential $\delta_E = \iota_\eta - \gamma \wedge \cdot$. This is a matrix factorization, and its endomorphism algebra in the category $\mathit{MF}(W)$ is the previously considered dga $\B$.

\begin{lemma}
In $D^b_{\mathrm{sing}}(W^{-1}(0))$, $\mathrm{coker}(\delta_E^1)$ is isomorphic to ${\EuScript S}_{W^{-1}(0),0}$.
\end{lemma}

\proof Take the chain complex of vector bundles on $W^{-1}(0)$ given by
\begin{equation} \label{eq:unrolled}
C^i = \begin{cases}
\bigoplus_{j = 0}^{[n/2]} \Omega^{2j-i}(V)|W^{-1}(0) & i \leq 0, \\
0 & i>0,
\end{cases}
\end{equation}
with $\delta_C^i = \iota_\eta - \gamma \wedge \cdot$ for $i<0$. This has the property that $\mathrm{coker}(\delta_E^1) = H^0(C)$. Consider the decreasing filtration $F_\bullet C$ whose pieces $F_m C$ consists of those summands in \eqref{eq:unrolled} with $j \geq m$. Passing to induced differential on the graded spaces of that filtration means that we keep only the $\iota_\eta$ term. In particular, $C/F_1C$ is the standard Koszul resolution of ${\EuScript S}_{V,0}$ restricted to $W^{-1}(0)$. Hence, its cohomology sheaves are the derived restrictions
\begin{equation}
H^i(C/F_1C) = {\EuScript S}_{V,0} \stackrel{\mathrm{L}^{-i}}{\otimes}_{\!\O(V)} \!\O(W^{-1}(0)) =
\begin{cases} {\EuScript S}_{W^{-1}(0),0} & i = 0,-1, \\ 0 & \text{otherwise.}\end{cases}
\end{equation}
The other quotients $F_mC/F_{m+1}C$, $m>0$, are truncations of the same resolution, which means that
\begin{equation} \label{eq:higher-coh}
H^i(F_mC/F_{m+1}C) = \begin{cases} \ker\big(\Omega^{2m}(V)|W^{-1}(0) \xrightarrow{\iota_\eta} \Omega^{2m-1}(V)|W^{-1}(0)\big) & i = 0,
\\ 0 &
i \neq 0. \end{cases}
\end{equation}
Therefore, $H^i(F_1C)$ is again zero in degrees $i \neq 0$, and $H^0(F_1C)$ is a successive extension of torsion-free sheaves, hence itself torsion-free. Finally, consider the long exact sequence
\begin{equation} \label{eq:long-ex}
\cdots H^{-1}(C/F_1C) \rightarrow H^0(F_1C) \rightarrow H^0(C) \rightarrow H^0(C/F_1C) \rightarrow 0.
\end{equation}
Because $H^{-1}(C/F_1C)$ is torsion, the leftmost arrow necessarily vanishes, which means that $H^0(C)$ is an extension of ${\EuScript S}_{W^{-1}(0),0}$ by $H^0(F_1C)$. We know that in the derived category $H^0(F_1C)$ is isomorphic to $F_1C$, which is a perfect complex, hence maps to zero when passing to $D^b_{\mathrm{sing}}(W^{-1}(0))$. This yields an isomorphism $H^0(C) \iso {\EuScript S}_{W^{-1}(0),0}$ in that category. \qed

We have now found a split-generator for the split-closure of $H^0(\mathit{MF}(W))$, and know that its endomorphism dga is $\B$. By the same general arguments as in Section \ref{sec:fukaya}, this implies:

\begin{corollary}
$D^\pi_{\mathrm{sing}}(W^{-1}(0)) \iso D^\pi(\B)$.
\end{corollary}

We will also need an equivariant version of this discussion, which is fairly straightforward. Suppose that $W$ and $\gamma$ are invariant under the action of a finite group $Z \subset \mathit{GL}(V)$. Defining equivariant categories $D^b_{sing,Z}(W^{-1}(0))$ and $\mathit{MF}_Z(W)$ in the obvious way, the analogue of Theorem \ref{th:orlov} holds \cite[Proposition 6.2]{quintero-velez08}. The skyscraper sheaf at the origin, with its natural equivariant structure, is no longer a split-generator. Instead, one should consider ${\EuScript S}_{W^{-1}(0),0} \otimes \C[Z]$, and the analogous equivariant matrix factorization $\Omega(V) \otimes \C[Z]$, whose endomorphism dga is the semidirect product $\B \rtimes Z$. The outcome is that
\begin{equation} \label{eq:equivariant-b}
D^\pi_{\mathrm{sing},Z}(W^{-1}(0)) \iso D^\pi(\B \rtimes Z).
\end{equation}

\section{The McKay correspondence\label{sec:mckay}}

We now return to $V = \C^3$ with the action of $\Z/5 \iso Z  \subset \mathit{SL}(V)$. The quotient has a canonical crepant resolution, namely the $G$-Hilbert scheme \cite{nakamura99}
\begin{equation} \label{eq:g-hilb}
X = \mathrm{Hilb}_Z(V) \longrightarrow \bar{X} = V/Z.
\end{equation}
The paper \cite{craw-reid99} gives an elementary toric description of $X$. Namely, take $N_\R = \R^3$, and let $N \subset N_\R$ be the lattice generated by $\Z^3$ together with $\frac{1}{5}{(1,1,3)}$. Let $\bar\Delta$ be the fan consisting of the single cone $N \cap \Z_+^3$ and its faces. This describes the affine toric variety $\bar{X}$. Now take the elementary simplex in $N_\R$, which is the one spanned by $\{(1,0,0),(0,1,0),(0,0,1)\}$, and triangulate it as follows:
\begin{equation}
\begin{centering}
\begin{picture}(0,0)%
\includegraphics{crepant.pstex}%
\end{picture}%
\setlength{\unitlength}{3947sp}%
\begingroup\makeatletter\ifx\SetFigFont\undefined%
\gdef\SetFigFont#1#2#3#4#5{%
  \reset@font\fontsize{#1}{#2pt}%
  \fontfamily{#3}\fontseries{#4}\fontshape{#5}%
  \selectfont}%
\fi\endgroup%
\begin{picture}(2636,2424)(1039,-1573)
\put(941,539){\makebox(0,0)[lb]{\smash{{\SetFigFont{10}{12.0}{\rmdefault}{\mddefault}{\updefault}{\color[rgb]{0,0,0}$(0,1,0)$}%
}}}}
\put(3076,-1336){\makebox(0,0)[lb]{\smash{{\SetFigFont{10}{12.0}{\rmdefault}{\mddefault}{\updefault}{\color[rgb]{0,0,0}$(1,0,0)$}%
}}}}
\put(941,-1336){\makebox(0,0)[lb]{\smash{{\SetFigFont{10}{12.0}{\rmdefault}{\mddefault}{\updefault}{\color[rgb]{0,0,0}$(0,0,1)$}%
}}}}
\end{picture}%
\end{centering}
\end{equation}
The cone over this triangulation yields a subdivision of $\bar\Delta$, whose fan $\Delta$ is the one describing $X$. Inspection of this picture shows that the preimage of the origin has two components, which are $\C P^2$ and a Hirzebruch surface $F_3$, intersecting each other in a rational curve.

We want to be pedestrian, and work through the standard construction of $X$ from $\Delta$. Let $M_\R = N_\R^\vee$ and $M = N^\vee$ be the dual space and dual lattice, respectively. Explicitly, $M_\R = \R^3$, and $M$ consists of those points $m \in \Z^3$ where $\frac{1}{5}(m_1 + m_2 + 3 m_3) \in \Z$. The cones $\sigma^\vee$ dual to the five maximal cones $\sigma \in \Delta$ are given by
\begin{equation}
\begin{aligned}
& \{m_2 \geq 0, \; m_1 + m_2 + 3 m_3 \geq 0, \; 2 m_1 + 2 m_2 + m_3 \geq 0\}, \\
& \{m_1 \geq 0, \; m_1 + m_2 + 3 m_3 \geq 0, \; 2 m_1 + 2 m_2 + m_3 \geq 0\}, \\
& \{m_2 \geq 0, \; m_3 \geq 0, \; m_1 + m_2 + 3 m_3 \geq 0\}, \\
& \{m_1 \geq 0, \; m_3 \geq 0, \; m_1 + m_2 + 3 m_3 \geq 0\}, \\
& \{m_1 \geq 0, \; m_2 \geq 0, \; 2 m_1 + 2 m_2 + m_3 \geq 0\}.
\end{aligned}
\end{equation}
When intersected with $M$, these give five semigroups $\sigma^\vee \cap M \iso \Z_+^3$, whose generators are
\begin{equation} \label{eq:generators}
\begin{aligned}
& \{(3,0,-1), \; (-1,0,2), \; (-1,1,0)\}, \\
& \{(0,3,-1), \; (0,-1,2), \; (1,-1,0)\}, \\
& \{(5,0,0), \; (-3,0,1), \; (-1,1,0)\}, \\
& \{(0,5,0), \; (0,-3,1), \; (1,-1,0)\}, \\
& \{(0,1,-2), \; (1,0,-2), \; (0,0,5)\}.
\end{aligned}
\end{equation}
These correspond to a covering of $X$ by toric charts, each of which is a copy of $\C^3$. Specifically, If one takes the generators listed in \eqref{eq:generators} to correspond to a basis of monomials in $\C[\sigma^\vee \cap M] \iso \C[a_1,a_2,a_3]$, the transformations from each chart to the previous one are
\begin{equation}
\begin{aligned}
& (a_1,a_2,a_3) \longmapsto (a_1a_3^3,a_2a_3^{-1},a_3^{-1}), \\
& (a_1,a_2,a_3) \longmapsto (a_2^{-1},a_1a_2^2,a_3), \\
& (a_1,a_2,a_3) \longmapsto (a_1a_3^5,a_2a_3^{-3},a_3^{-1}), \\
& (a_1,a_2,a_3) \longmapsto (a_1^5 a_3^2, a_1^{-3} a_3^{-1}, a_1^{-1} a_2).
\end{aligned}
\end{equation}

The $Z$-invariant hypersurface $W^{-1}(0)$ descends to $\bar{H} \subset \bar{X}$, which is a singular surface that turns out to be rational.
%: an explicit birational map $\bar{H} \rightarrow \C^2$ is $x \mapsto (x_1^{-1} x_3^2, x_1^{-1}x_2)$.
We are interested in its preimage under the resolution \eqref{eq:g-hilb}, denoted by $H$. In the charts constructed above, the defining equation for $H$ is, respectively,
\begin{equation}
\begin{aligned}
& a_1a_2(a_3 - a_1 - a_1 a_3^5 - a_2^2) = 0, \\
& a_1a_2(a_3 - a_1 a_3^5 - a_1 - a_2^2) = 0, \\
& a_1(a_2a_3 - 1 - a_3^5 - a_1^2 a_2^5) = 0, \\
& a_1(a_2a_3 - a_3^5 - 1 - a_1^2 a_2^5) = 0, \\
& a_3(a_1a_2 - a_2^5 a_3 - a_1^5 a_3 - 1) = 0.
\end{aligned}
\end{equation}
This is a surface with three components $H_1,H_2,H_3$ and only normal crossing singularities. The first component is $H_1 \iso \C P^2$, which is $\{a_1 = 0\}$ in the first and second charts, and $\{a_3 = 0\}$ in the fifth chart. The intersections $H_{12},H_{13} \subset H_1$ are a line and smooth conic, respectively, which are in general position. The second component is $H_2 \iso F_3$, which is $\{a_2 = 0\}$ in the first and second charts, and $\{a_1 =0\}$ in the third and fourth charts. Let's explicitly identify $H_2$ with the ruled surface $P(\O \oplus \O(-3)) \rightarrow \C P^1$, where we have sections $S_- = P(\O \oplus \{0\})$ and $S_+ = P(\{0\} \oplus \O(-3))$ with selfintersection $\pm 3$. Then $H_{12} = S_-$, while $H_{23}$ is a section intersecting $S_-$ transversally over $0, \infty \in \C P ^1$, and intersecting $S_+$ transversally over the five points $[1:\eta] \in \C P^1$, $\eta \in \sqrt[5]{1}$. These properties determine $H_{23} \subset H_2$ uniquely up to the fibrewise $\C^*$-action on the ruling. In particular, it has selfintersection number $7$. The last component is the non-toric one, which is the proper transform of $\bar{H}$. One can analyze this through the map $H_3 \rightarrow \C P^2$, given by $(a_1,a_2,a_3) \mapsto [1:i a_2:-a_3]$ in the first of our five charts. Inside $\C P^2$ take again a smooth conic and a line, which concretely are given by $z_1^2 + z_0z_2 = 0$ and $z_1 = 0$ in homogeneous coordinates $[z_0:z_1:z_2]$. Now blow up the ten points $[1:\zeta:\zeta^2]$ for $\zeta \in \sqrt[10]{1}$. Furthermore, blow up the five points $[1:0:\eta]$ for $\eta \in \sqrt[5]{1}$, and another five points infinitely close to them, corresponding to functions whose derivative vanishes along our line. The outcome of this blowup process is a compactification of $H_3$. The intersections $H_{13}$ and $H_{23}$ are the proper transforms of our conic and line, respectively. In particular, their selfintersection numbers are $-6$ and $-9$, respectively.

Consider the category $D^b_{\mathrm{sing}}(H)$ of Landau-Ginzburg branes, defined as in \eqref{eq:db-sing}, and its split-closure $D^\pi_{\mathrm{sing}}(H)$.

\begin{theorem} \label{th:derived-mckay}
There is an equivalence of triangulated categories,
\begin{equation}
D^\pi_{\mathrm{sing},Z}(W^{-1}(0)) \iso D^\pi_{\mathrm{sing}}(H).
\end{equation}
\end{theorem}

The existence of a full and faithful functor $D^b_{\mathrm{sing},Z}(W^{-1}(0)) \rightarrow D^\pi_{\mathrm{sing}}(H)$ was proved in \cite{mehrotra05}. Essential surjectivity is proved in \cite{quintero-velez08} for a somewhat different special case, but the method given there adapts to our situation. In fact, since we only need the statement for the split-closures, it is enough for us to consider complexes of sheaves with compactly supported cohomology, which removes the need for the most technical aspects of \cite{bridgeland-king-reid01}. From \eqref{eq:equivariant-b} and Proposition \ref{th:endomorphism-characterized}, it follows that $D^\pi_{\mathrm{sing},Z}(W^{-1}(0)) \iso D^\pi(\A \rtimes Z)$, where $\A$ is as in Proposition \ref{th:classification}. The conclusion is:

\begin{corollary}
$D^\pi_{\mathrm{sing}}(H) \iso D^\pi(\A \rtimes Z)$.
\end{corollary}

The right hand side is the same as in the description of the Fukaya category given in Corollary \ref{th:describes-f}. Hence, the combination of the two results implies Theorem \ref{th:main}.


\begin{thebibliography}{10}

\bibitem{abouzaid08}
M.~Abouzaid.
\newblock On the {F}ukaya categories of higher genus surfaces.
\newblock {\em Adv. Math.}, 217:1192--1235, 2008.

\bibitem{akaho05}
M.~Akaho.
\newblock Intersection theory for {L}agrangian immersions.
\newblock {\em Math. Res. Lett.}, 12:543--550, 2005.

\bibitem{akaho-joyce08}
M.~Akaho and D.~Joyce.
\newblock Immersed {L}agrangian {F}loer theory.
\newblock Preprint arXiv:0803.0717.

\bibitem{arnold-gusein-zade-varchenko}
V.~I. Arnold, S.~M. Gusein-Zade, and A.~N. Varchenko.
\newblock {\em Singularities of differentiable maps}.
\newblock Birkh{\"a}user, 1988.

\bibitem{aspinwall-katz06}
P.~Aspinwall and S~Katz.
\newblock Computation of superpotentials for {D}-branes.
\newblock {\em Comm. Math. Phys.}, 264:227--253, 2006.

\bibitem{atiyah57}
M.~Atiyah.
\newblock Vector bundles over an elliptic curve.
\newblock {\em Proc. London Math. Soc.} 7:414--452, 1957.

\bibitem{auroux-katzarkov-orlov04}
D.~Auroux, L.~Katzarkov and D.~Orlov.
\newblock Mirror symmetry for weighted projective planes and their noncommutative deformations.
\newblock {\em Ann. of Math. (2)}, 167:867--943, 2008.

\bibitem{auroux07}
D.~Auroux.
\newblock Mirror symmetry and {T}-duality in the complement of an anticanonical
  divisor.
\newblock Preprint arXiv:0706.3207.

\bibitem{balmer-schlichting01}
P.~Balmer and M.~Schlichting.
\newblock Idempotent completion of triangulated categories.
\newblock {\em J. Algebra}, 236:819--834, 2001.

\bibitem{baranovsky-pecharich09}
V.~Baranovsky and J.~Pecharich.
\newblock On equivalences of derived and singular categories.
\newblock Preprint arXiv:0907.3717, 2009.

\bibitem{beilinson-ginzburg-soergel96}
A.~Beilinson, V.~Ginzburg, and W.~Soergel.
\newblock Koszul duality patterns in representation theory.
\newblock {\em Journal of the Amer. Math. Soc.}, 9:473--527, 1996.

\bibitem{birman72}
J.~Birman.
\newblock A normal form in the homeotopy group of a surface of genus 2, with
  applications to 3-manifolds.
\newblock {\em Proceedings of the Amer. Math. Soc.}, 34:379--384, 1972.

\bibitem{bridgeland-king-reid01}
T.~Bridgeland, A.~King, and M.~Reid.
\newblock The {M}c{K}ay correspondence as an equivalence of derived categories.
\newblock {\em Journal of the Amer. Math. Soc.}, 14:535--554, 2001.

\bibitem{chekanov99}
Yu.~Chekanov.
\newblock Differential algebras of Legendrian links.
\newblock {\em Invent. Math.}, 150:441--483, 2002.

\bibitem{cho02}
C.-H. Cho.
\newblock Holomorphic discs, spin structures, and the {F}loer cohomology of the
  {C}lifford torus.
\newblock {\em Intern. Math. Res. Notices}, 35:1803--1843, 2004.

\bibitem{cho-oh02}
C.-H. Cho and Y.-G. Oh.
\newblock Floer cohomology and disc instantons of {L}agrangian torus fibers in
  toric {F}ano manifolds.
\newblock {\em Asian J. Math.}, 10:773--814, 2006.

\bibitem{cornea-lalonde06}
O.~Cornea and F.~Lalonde.
\newblock Cluster homology: an overview of the construction and results.
\newblock {\em Electron. Res. Announc. Amer. Math. Soc.}, 12:1--12
  (electronic), 2006.

\bibitem{craw-reid99}
A.~Craw and M.~Reid.
\newblock How to calculate {$A$}-{H}ilb {$\mathbb{C}^3$}.
\newblock In {\em Geometry of toric varieties}, volume~6 of {\em S\'emin.
  Congr.}, pages 129--154. Soc. Math. France, 2002.

\bibitem{dyckerhoff09}
T.~Dyckerhoff.
\newblock Compact generators in categories of matrix factorizations.
\newblock Preprint arXiv:0904.4713.

\bibitem{efimov09}
A.~Efimov.
\newblock Homological mirror symmetry for curves of higher genus.
\newblock Preprint arXiv:0907.3903, 2009.

\bibitem{eisenbud80}
D.~Eisenbud.
\newblock Homological algebra on a complete intersection, with an application
  to group representations.
\newblock {\em Trans. Amer. Math. Soc.}, 260:35--64, 1980.

%\bibitem{fukaya-oh98}
%K.~Fukaya and Y.-G. Oh.
%\newblock Zero-loop open strings in the cotangent bundle and {M}orse homotopy.
%\newblock {\em Asian J. Math.}, 1:96--180, 1998.
%
\bibitem{fooo08}
K.~Fukaya, Y.-G. Oh, H.~Ohta, and K.~Ono.
\newblock Lagrangian {F}loer theory on compact toric manifolds. {I}.
\newblock Preprint arXiv:0802.1703.

\bibitem{fooo08b}
K.~Fukaya, Y.-G. Oh, H.~Ohta, and K.~Ono.
\newblock Lagrangian {F}loer theory on compact toric manifolds. {II}: {B}ulk
  deformations.
\newblock Preprint arXiv:0810.5654.

\bibitem{fooo}
K.~Fukaya, Y.-G. Oh, H.~Ohta, and K.~Ono.
\newblock Lagrangian intersection {F}loer theory - anomaly and obstruction.
\newblock In press.

\bibitem{fukaya02b}
Kenji Fukaya.
\newblock Mirror symmetry of abelian varieties and multi-theta functions.
\newblock {\em J. Algebraic Geom.}, 11:393--512, 2002.

\bibitem{4}
B.~Fang, C.-C. Liu, D. Treumann and E. Zaslow.
\newblock The coherent-constructible correspondence and Homological Mirror Symmetry for toric varieties.
\newblock Preprint arXiv:0901.4276, 2009.

\bibitem{goldman-millson88}
W.~Goldman and J.~Millson.
\newblock The deformation theory of the fundamental group of compact
  {K}{\"a}hler manifolds.
\newblock {\em IHES Publ. Math.}, 67:43--96, 1988.

\bibitem{hori-vafa00}
K.~Hori and C.~Vafa.
\newblock Mirror symmetry.
\newblock Preprint hep-th/0002222, 2000.

\bibitem{horja}
P.~Horja.
\newblock Derived category automorphisms from mirror symmetry.
\newblock {\em Duke Math. J.} 127:1--34, 2005.

\bibitem{hochschild-kostant-rosenberg62}
G.~Hochschild, B.~Kostant, and A.~Rosenberg.
\newblock Differential forms on regular affine algebras.
\newblock {\em Trans. Amer. Math. Soc.}, 102:383--408, 1962.

\bibitem{huybrechts-thomas}
D.~Huybrechts and R.~Thomas.
\newblock $\mathbb{P}$-objects and autoequivalences of derived categories.
\newblock {\em Math. Res. Lett.} 13:87--98, 2006.

%\bibitem{ito-nakamura95}
%Y.~Ito and I.~Nakamura.
%\newblock Mc{K}ay correspondence and {H}ilbert schemes.
%\newblock {\em Proc. Japan Acad. Ser. A Math. Sci.}, 72:135--138, 1996.
%
%\bibitem{joyce08}
%D.~Joyce.
%\newblock {K}uranishi homology and {K}uranishi cohomology.
%\newblock Preprint arXiv:0707.3572.

\bibitem{katzarkov07}
L.~Katzarkov.
\newblock Birational geometry and homological mirror symmetry.
\newblock In {\em Real and complex singularities}, pages 176--206. World Sci.
  Publ., 2007.

\bibitem{kapustin-katzarkov-orlov-yotov09}
A.~Kapustin, L.~Katzarkov, D.~Olov and M.~Yotov.
\newblock Homological Mirror Symmetry for manifolds of general type.
\newblock Preprint, 2009.

\bibitem{kontsevich97}
M.~Kontsevich.
\newblock Deformation quantization of {P}oisson manifolds {I}.
\newblock {\em Lett. Math. Phys.}, 66:157--216, 2003.

\bibitem{kontsevich94}
M.~Kontsevich.
\newblock Homological algebra of mirror symmetry.
\newblock In {\em Proceedings of the International Congress of Mathematicians
  (Z{\"u}rich, 1994)}, pages 120--139. Birkh{\"a}user, 1995.

\bibitem{kontsevich-soibelman00}
M.~Kontsevich and Y.~Soibelman.
\newblock Homological mirror symmetry and torus fibrations.
\newblock In {\em Symplectic geometry and mirror symmetry}, pages 203--263.
  World Scientific, 2001.

\bibitem{katzarkov-kontsevich-pantev08}
T.~Pantev, L.~Katzarkov, M.~Kontsevich.
\newblock Hodge theoretic aspects of mirror symmetry.
\newblock Preprint arXiv:0806.0107.

\bibitem{keller99}
B.~Keller.
\newblock Introduction to {$A$}-infinity algebras and modules.
\newblock {\em Homology Homotopy Appl. (electronic)},
3:1--35, 2001.

\bibitem{lada-markl95}
T.~Lada and M.~Markl.
\newblock Strongly homotopy {L}ie algebras.
\newblock {\em Comm. Algebra}, 23:2147--2161, 1995.

\bibitem{markl04}
M.~Markl.
\newblock Transferring {$A_\infty$} (strongly homotopy associative) structures.
\newblock {\em Rend. Circ. Mat. Palermo (2) Suppl.}, 79:139--151, 2006.

\bibitem{mcduff90b}
D.~McDuff.
\newblock The structure of rational and ruled symplectic $4$-manifolds.
\newblock {\em Journ. Amer. Math. Soc.}, 3:679--712, 1990.
\newblock Erratum: same journal, 5:987--988, 1992.

\bibitem{mehrotra05}
S.~Mehrotra.
\newblock {\em Triangulated categories of singularities, matrix factorizations
  and {LG} models}.
\newblock PhD thesis, Univ. of Pennsylvania, 2005.

\bibitem{mukai}
S.~Mukai.
\newblock Duality between $D(X)$ and $D(\hat X)$ with its application to Picard sheaves.
\newblock {\em Nagoya Math. J.} 81:153--175, 1981.

\bibitem{nakamura99}
I.~Nakamura.
\newblock Hilbert schemes of abelian group orbits.
\newblock {\em J. Algebraic Geom.}, 10:757--779, 2001.

\bibitem{oh93}
Y.-G. Oh.
\newblock Floer cohomology of {L}agrangian intersections and pseudo-holomorphic
  discs {I}.
\newblock {\em Comm. Pure Appl. Math.}, 46:949--994, 1993.

\bibitem{orlov08}
D.~Orlov.
\newblock {D}-branes of type {B} in {L}andau-{G}inzburg models and mirror
  symmetry for the genus two curve.
\newblock Conference talk, IAS Princeton, 3/12/2008.

\bibitem{orlov08b}
D.~Orlov.
\newblock {F}ormal completions and idempotent completions of
triangulated categories of singularities.
\newblock Preprint arXiv:0901.1859.

\bibitem{orlov04}
D.~Orlov.
\newblock Triangulated categories of singularities and {D}-branes in
  {L}andau-{G}inzburg models.
\newblock {\em Proc. Steklov Inst. Math.}, 246:227--248, 2004.

\bibitem{polishchuk-zaslow98}
A.~Polishchuk and E.~Zaslow.
\newblock Categorical mirror symmetry: the elliptic curve.
\newblock {\em Adv. Theor. Math. Phys.}, 2:443--470, 1998.

\bibitem{polishchuk98}
A.~Polishchuk.
\newblock Massey and Fukaya products on elliptic curves.
\newblock {\em Adv. Theor. Math. Phys.}, 4:1187--1207,2000.

\bibitem{rossi08}
P.~Rossi.
\newblock {G}romov-{W}itten theory of orbicurves, the space of tri-polynomials
  and {S}ymplectic {F}ield {T}heory of {S}eifert fibrations.
\newblock Preprint arXiv:0808.2626.

\bibitem{schwarz96}
M.~Schwarz.
\newblock A quantum cup-length estimate for symplectic fixed points.
\newblock {\em Invent. Math.}, 133:353--397, 1998.

\bibitem{seidel08x}
P.~Seidel.
\newblock Python code for computing $A_\infty$-deformations, 2008. Can be
  downloaded from the author's homepage, {\tt
  http://math.mit.edu/{\~{}}seidel}.

\bibitem{seidel03b}
P.~Seidel.
\newblock Homological mirror symmetry for the quartic surface. \newline
\newblock Preprint math.SG/0310414, 2003.

\bibitem{seidel01}
P.~Seidel.
\newblock A long exact sequence for symplectic {F}loer cohomology.
\newblock {\em Topology}, 42:1003--1063, 2003.

\bibitem{seidel04}
P.~Seidel.
\newblock {\em {F}ukaya categories and {P}icard-{L}efschetz theory}.
\newblock European Math. Soc., 2008.

\bibitem{seidel-thomas}
P.~Seidel and R.~Thomas.
\newblock {\em {B}raid group actions on derived categories of coherent sheaves}.
\newblock {\em Duke Math. Jour.}, 108:37--108, 2001.

\bibitem{subotic}
A.~Subotic.
\newblock Ph.D. thesis, MIT, in preparation.

\bibitem{takahashi07}
A.~Takahashi.
\newblock Weighted projective lines associated to regular systems of weights of dual type.
\newblock Preprint arXiv:0711.3907.

\bibitem{tougeron70}
J.-C. Tougeron.
\newblock Id\'eaux des fonctions diff\'erentiables.
\newblock {\em Ann. Inst. Fourier}, 18:177--240, 1970.

\bibitem{quintero-velez08}
A.~Quintero Velez.
\newblock Mc{K}ay correspondence for {L}andau-{G}inzburg models.
\newblock Preprint arXiv:0711.4774.

\bibitem{wehrheim-woodward06}
K.~Wehrheim and C.~Woodward.
\newblock Functoriality for {L}agrangian correspondences in {F}loer theory.
\newblock Preprint arXiv:0708.2851.

\bibitem{weibel}
C.~A. Weibel.
\newblock {\em An introduction to homological algebra}.
\newblock Cambridge Univ. Press, 1994.

\end{thebibliography}
\end{document}